\documentclass[12pt]{amsart}
\usepackage[applemac]{inputenc}
\usepackage{amsmath,amsthm,amssymb,latexsym,wasysym}
\usepackage[T1]{fontenc}
\usepackage{libertine} % para texto
\usepackage{euler} %para matem\'aticas
\usepackage{eucal} %para matem\'atica - caligr.
\usepackage[]{color}
\usepackage[]{graphicx}
\usepackage{pb-diagram}
\usepackage{pst-3d,pst-coil,pst-eps,pst-fill,pst-grad,pst-math,pst-text,pst-tree}
\newcommand\C[1]{\mathcal{#1}}
\newcommand\M[1]{\mathfrak{#1}}

\def\monster{\mathbb{M}}
\def\E{\varepsilon}
\def\dist{\mathbf{d}}
\newcommand{\indep}{\makebox[13pt]{\raisebox{-0.75ex}{$\smile$}\hspace{-3.00mm}
$\mid$\, }}
\newcommand{\eindep}{\makebox[13pt]{\raisebox{-0.75ex}{$\smile$}\hspace{-3.00mm}
$\mid$\, }^{\hspace{-1mm}\varepsilon}}

\def\ordencea{\prec_{\C{K}}}
\def\ordenceac{\succ_{\C{K}}}
\newcommand{\eop}[1]{
\hspace{10mm} \vspace{-6mm}
\begin{flushright}
\qedsymbol$_{\text{#1}}$\\ \ \\
\end{flushright}
}
\newenvironment{prueba}[1][{\it Proof}]{\noindent {\it #1.} }{}
\newcommand{\bdem}[1][Proof]{\begin{prueba}[#1]}
\newcommand{\edem}[1][]{\eop{#1}
\end{prueba}}
\def\bsdem{\begin{prueba}[Reference]}\def\bsindem{\begin{proof}[\ ]}
\def\tuple{\overline}
\def\rest{\upharpoonright}

\mathchardef\mhyphen="2D
\newcommand{\gatp}{{ga\mhyphen tp}}
\newcommand{\gaS}{{ga\mhyphen S}}

\newcommand{\stype}{\mathfrak{S}t}
\newcommand{\tower}{(\mathfrak{M},\overline{a},\mathfrak{N},\C{M}, \C{N})}

\newtheorem{theorem}{Theorem}[subsection] 
\newtheorem{lemma}[theorem]{Lemma} 
\newtheorem{proposition}[theorem]{Proposition} 
\newtheorem{corollary}[theorem]{Corollary} 
 
\newtheorem{definition}[theorem]{Definition}

\newtheorem{remark}[theorem]{Remark}

\newtheorem{fact}[theorem]{Fact} 
\newtheorem{assumption}[theorem]{Assumption} 
\newtheorem{notation}[theorem]{Notation}

\title{Limit Models in Metric Abstract Elementary Classes: the
  categorical case}

\author{Andr\'es Villaveces and Pedro Zambrano\\
Universidad Nacional de Colombia - Bogotá}
\begin{document}

\thanks{The second author wants
  to thank the first author for the time devoted to advising him for
  his Ph.D.'s thesis, one of whose fruits is this paper. Both authors
  were partially supported by Colciencias grant \emph{M\'etodos de
    estabilidad en clases no estables.}\\
e-mail: {\sf
       avillavecesn@unal.edu.co}, {e-mail: {\sf
       phzambranor@unal.edu.co} - Departamento de Matem\'aticas
     Universidad Nacional de Colombia, AK 30 $\#$ 45-03, Bogot\'a
     111321 - Colombia}}
% \author[P. Zambrano]{Pedro Zambrano\footnote{e-mail: {\sf phzambranor@unal.edu.co} \\Departamento de Matem\'aticas Universidad Nacional de Colombia, AK 30 $\#$ 45-03 c\'odigo postal 111321, Bogot\'a - Colombia}}

\maketitle

\begin{abstract}
We study versions of limit models adapted to the context of
\emph{metric abstract  elementary classes}. Under categoricity and
superstability-like assumptions, we generalize some
theorems from {\cite{GrVaVi,ShVi,Va,Va2}}. We prove
criteria for existence and uniqueness of
limit models in the me\-tric context.
\end{abstract}

\section{Preliminaries - Why metric AECs, limit models and towers?}

In his seminal paper published in 1975~\cite{ShLazyGuide}, Shelah
sketched four ``kinds of classes for which we have stability theory'' -
four early extensions of First
Order logic, in the more general context of the existence of large
homogeneous models. One of the four kinds (Kind III - existentially
closed models of a universal theory with JEP) included a particular
case (stability theory for normed spaces and for Banach spaces, only
hinted at in~\cite{ShLazyGuide}). A long array
of results, starting with earlier work by Chang and
Keisler~\cite{ChKe}, with crucial ideas from Henson and
Iovino~\cite{HeIo}, finally converged in the
monograph~\cite{CoMon} that settled the foundation of First Order
metric continuous model theory. Our work aims at wider classes of
models, and builds on a combination of Abstract Elementary Classes
with more purely metric Continuous Logic constructions. It is a
contribution toward the stability theory of Metric Abstract Elementary
Classes. 

The Model Theory of metric structures can be generalized %marzo de 2013
 in an effec\-tive %marzo de 2013 
way to the Abstract Elementary Class (AEC) context %marzo de 2013
by blending some of the constructions typical of AECs with ideas and
paradigms from First Order
Continuous Model Theory as understood by
%extending in some senses
%the framework of First Order
%Continuous Model Theory
~\cite{CoMon} and in other senses benefitting from the enormous wealth
% marzo de 2013
 of Stability Theory in Abstract Elementary Classes. Other
authors (Hirvonen~\cite{Hi} in her thesis with Hyttinen, and also
independently Shelah and Usvyatsov~\cite{ShUs}) have provided other
frameworks for model theory of metric structures outside continuous
first order.

Hirvonen and Hyttinen have developed a solid framework for
ca\-te\-go\-ri\-ci\-ty transfer of metric AEC and for the study of
$\aleph_0$-stable classes of metric structures (a good analysis of
primary models, basic items in the definition itself, etc.).
%They provide a good analysis of ... llenar

Our focus here is the beginning of %marzo 2013
an analysis of %``superstability'' PZ: junio 7, 2014.
weak notion of superstability
%\marginpar{\textcolor{red}{PZ: jun 7, 2014.}}
in metric AEC. Of course this goal is long-winded, but we provide first
steps in that direction in this paper. In particular, building mainly
on ideas from the discrete AEC setting coming from~{\cite{GrVaVi,
    ShVi, Va,Va2},} and
related more distantly to Shelah's ideas in~\cite{Sh600}, we approach
here the connection between two facets of (protean) superstability:
limit models (existence and first steps towards uniqueness).

%\marginpar{\textcolor{red}{PZ:jun7}}
{One of the main aims of the study of uniqueness of
  limit models in \cite{GrVaVi,ShVi,Va} is proving a categoricity
  transfer theorem in their settings. Also, J. Baldwin~\cite{Ba2} has
  an unpublished work about a rough notion of domination of Galois
  types in AECs under uniqueness of limit models. Although Hirvonen
  studied categoricity in $\omega$-stable homogeneous MAECs~\cite{Hi},
  we do not focus on such kind of results in our work. Our interest to
  study uniqueness of limit models is understanding this notion as a
  weak version of superstability in MAECs and some of its consequences
  in domination, orthogonality and parallelism of Galois types in
  MAECs~\cite{ViZa1}.
}

The main constructions in our paper are versions of towers adapted to
the metric context (s-towers and metric s-towers). Specifically,
reduced towers have been
used extensively by  Shelah and Villaveces~\cite{ShVi}, Grossberg,
VanDieren and Villaveces~\cite{GrVaVi}, and Jarden and
Shelah~\cite{JaSh} in their work on
AEC before. Here we adapt them to the metric setting and use them to
prove various lemmas useful to an approach of
uniqueness of limit models in metric AEC (in a forthcoming
paper). Towers can be regarded as a strong generalization of the
concept of Galois type: a Galois type is (an equivalence class) of
triples $(M,N,a)$ where $M\prec N$ and $a\in N\setminus M$ --- towers
``refine'' the way the element $a$ is connected to $M$ inside $N$ and
provides a very robust situation where $a$ is replaced by a long
sequence $(a_i)_i$, and the models $M$ and $N$ themselves are ``sliced
through''. Some of our results on towers may be regarded as the
analogue of extension properties of Galois types: the beginning of a
generalization of
 independence and ``forking'' calculus-like properties of the triples
 underlying Galois types.
% may be lifted in a robust way to  towers. This has been explored by
% the authors mentioned above in the usual AEC setting - we begin the
% exploration here of the \emph{metric} version.
%Agregar refs

Towers play a key role in the proof of uniqueness of limit models%\marginpar{\textcolor{red}{PZ:jun8}}
given in \cite{GrVaVi} as they allow us to build models that are
simultaneously $\theta$-limits and $\omega$-limits: we will construct a
$(\mu,\theta)$-tower which is also a $(\mu,\omega)$-tower for every
limit ordinal $\theta<\mu^+$ by building a \emph{rectangular} array
witnessing these two ways of being a limit, for the same model,
ultimately. In this paper, we prove -under
$\mu^+$-categoricity and a weak version of superstability-  {\it
  uniqueness -up to isomorphism- of limit
  models of density character $\mu$} (i.e., if $M_1$ is a
$(\mu,\theta_1)$-limit model over $M$ and $M_2$ is a
$(\mu,\theta_2)$-limit model over $M$ with $\theta_1,\theta_2<\mu^+$
limit ordinals and $dc(M_1)=dc(M_2)$, then
$M_1\approx_M M_2$) under suitable superstability-like assumptions. If
$cf(\theta_1)=cf(\theta_2)$, then by a standard {\it back and forth}
argument we are done.  So, if $cf(\theta_1)\neq cf(\theta_2)$, as in
\cite{GrVaVi}, the key idea is to build a {$(\mu,\theta_i)$-limit model
over $M$ $M_{\theta_i}$ ($i\in \{1,2\}$)}
%\marginpar{\textcolor{red}{PZ: jun7}}
which is also a $(\mu,\omega)$-limit model over
$M$ for any ordinal $\theta<\mu^+$, so 
\begin{eqnarray*}
M_1 &\approx_M& M_{\theta_1} \text{\ (because they are
  $(\mu,\theta_1)$-limits over $M$)}\\\
&\approx_M & M_{\theta_2} \text{\ (because they are
  $(\mu,\omega)$-limits over $M$)}\\
&\approx_M &M_2 \text{\ (because they are $(\mu,\theta_2)$-limits over
  $M$)}
\end{eqnarray*}

%In order to build a model such that $(\mu,\theta)$-limit model over $M$ $M_\theta$ which is also a $(\mu,\omega)$-limit model over $M$ for any ordinal $\theta<mu^+$,
%PZ: jun 7 2014
%\marginpar{\textcolor{red}{PZ:jun7}}
Given any ordinal $\theta<\mu^+$,  as in \cite{GrVaVi}, to be able to construct a $(\mu,\theta)$-limit model over $M$ $M_\theta$ which is also a $(\mu,\omega)$-limit model over $M$,  we define the notion
of {\it smooth tower}, which corresponds to an adaptation of the
notion of {\it tower} given in \cite{GrVaVi} but in our metric
setting. The key idea is to extend (via $\C{K}$-embeddings) a given
tower of length of cofinality $\theta$ to a  special kind of tower
({\it reduced towers}) which is continuous and to a kind of tower
({\it relatively full tower}) which satisfies a kind of relative
saturation. Iterating this argument $\omega$ many times, the idea is
to prove that the directed limit of such directed system is  a {\it
  reduced} (and therefore a continuous) tower where the completion of
its union is a $(\mu,\theta)$-limit model over $M$ (which is
consequence of the {\it full-relativeness} of the extensions given in
the directed system). To be a $(\mu,\omega)$-limit model over $M$ is
assured defining in a suitable way the notion of extension of
towers. Although this argument has the same general outline as the
proof done
in \cite{GrVaVi}, we point out that many details in our proof here are
quite
different: e.g., our notion of {\it reduced towers} involves a
Lipschitz-like function which determines a notion of closeness  of
towers instead of intersections as in \cite{GrVaVi} and we have to
deal with completion of union of increasing chains of towers in the
metric sense instead of just its union, which makes more complicated
some of the arguments if we compare them with the proofs given in
\cite{GrVaVi}.

In \cite{GrVaVi} the authors prove the uniqueness of limit models
under superstability-like assumptions for AEC. Here we study the
behavior of s-towers under superstability-like assumptions
(assumption~\ref{Superstability}) and categoricity for the metric
setting. Remark~\ref{rem-categoricity} at the beginning of the proof
and the use of the main hypothesis addresses in detail some of the
historic issues with the use of categoricity in proofs of uniqueness
of limit models.

%%completar

\emph{Note}: we base many of our results here in the constructions
of~\cite{ViZa1}, where we define and study the notions of
$\varepsilon$-splitting and smooth independence -- both of them
notions of independence appropriate for the metric AEC case. We study
there conditions for good
behavior of stationarity, existence, extension, etc. We use freely
these notions in this paper, and refer to their statement
in~\cite{ViZa1} at crucial places here.

We would like to thank John Baldwin, Rami Grossberg, Tapani Hyttinen
and Monica VanDieren for extremely helpful discussions on our proof
and on the various concepts presented here. We also want to thank the
referee for many questions and suggestions that have led to serious
improvements and clarifications in this paper.

\section{Basic facts on Metric AECs}

We consider a natural adaptation of the notion of {\it Abstract
Elementary Class} (see \cite{Gr} and \cite{BaMon}), but work in
a context of Continuous Logic that generalizes the ``First Order
Continuous'' setting of~\cite{CoMon} by removing the assumption of
uniform continuity\footnote{Uniform continuity guarantees logical
  compactness in their formalization, but we drop compactness in
  AEC-like settings.}. We base our definitions on \cite{Hi,GrVa}.
  
\begin{definition}\label{DensityCharacter}
The \emph{density character} of a topological space is the smallest
cardinality of a dense subset of the space. If $X$ is a topological
space, we denote its density character by $dc(X)$. If $A$ is a
subset of a topological space $X$, we define
$dc(A):=dc(\overline{A})$.
%:
\end{definition}

%\marginpar{\textcolor{red}{PZ:jun8}. Lo coloqué porque a pesar de que es una definición estándar, el referee estaba confundido porque aparecía en varias partes del artículo y dijo que no se había definido. Entonces creo que mejor hacerlo.}
\begin{definition}\label{dist_element_set}
Let $(M,d)$ be a metric space, $X\subseteq M$ and $a\in M$. The \emph{distance} between
$a$ and $X$ is defined as $d(a,X):=\inf\{d(a,b): b\in X\}$.
\end{definition}

\begin{definition}\label{MAEC}\label{maec}
Let $\mathcal{K}$ be a class of $L$-structures (as in first order
continuous logic, all the metric spaces $(M,d)$ considered here are
complete in the
metric sense; the interpretations of function symbols $F^n$ of arity
$n<\omega$ are continuous - here not necessarily uniformly continuous
as in first order continuous logic, as we do not strive for compactness;
% functions of the form $F^M:M^n\to M$,
the interpretations of relational symbols $R^k$ of arity $k<\omega$ are
continuous functions $R^M: M^k\to [0,1]$ and the interpretations of
constant symbols are elements in $M$) and let $\ordencea$ be a binary
relation defined in
$\C{K}$. We say that $(\C{K},\ordencea)$ is a {\it Metric Abstract
Elementary Class} (shortly {\it MAEC}) if:

\begin{enumerate}
\item $\C{K}$ and $\ordencea$ are closed under isomorphism.% $\cong$.
\item $\ordencea$ is a partial order in $\C{K}$.
\item If $M\ordencea N$ then $M\subseteq N$.
\item (\emph{Completion of Union of Chains}) If $(M_i:i<\lambda)$ is a $\ordencea$-increasing chain then
    \begin{enumerate}
    \item the function symbols in $L$ can be uniquely interpreted on the completion of
        $\bigcup_{i<\lambda} M_i$ in such a way that
        $\overline{\bigcup_{i<\lambda}M_i} \in \C{K}$
    \item for each $j < \lambda$ , $M_j \ordencea
        \overline{\bigcup_{i<\lambda} M_i}$
    \item if for each $i<\lambda$ $M_i\ordencea N$, then
        $\overline{\bigcup_{i<\lambda} M_i} \ordencea N$.
    \end{enumerate}
    \item (\emph{Coherence}) if $M_1\subseteq M_2\ordencea M_3$ and
        $M_1\ordencea M_3$, then $M_1\ordencea M_2$.
    \item (DLS) There exists a cardinality $LS(K)$ (which is called the
      {\it metric L\"o\-wen\-heim-Skolem number}) such that if $M
      \in \C{K}$
        and $A \subseteq M$, then there exists $N\in \C{K}$ such that $dc(N) \le dc(A) + LS(K)$
        and $A\subseteq N \ordencea M$.
\end{enumerate}
\end{definition}

In Definition~\ref{maec} (4) we take
  \emph{completions} of unions at limit stages
  (see~\cite{CoMon}); unions of $\omega$-chains of complete metric
  spaces are not necessarily complete, hence our requirement.
% if we take a union of an increasing chain of
%   Banach spaces of  length $\omega$ $\langle M_n:n<\omega \rangle$, it
%   is possible to find a cofinal Cauchy sequence in $\bigcup_{n<\omega}
%   M_n$ (i.e., a Cauchy sequence $(a_n)_{n<\omega}$ in
%   $\bigcup_{n<\omega} M_n$ such that for every $k<\omega$ there exists
%   $m<\omega$ such that $a_m\in M_k$) which does not converges in
%   $\bigcup_{n<\omega} M_n$. For being able to obtain a complete metric
%   space, we have to take its closure $\overline{\bigcup_{n<\omega}
%     M_n}$ instead of just $\bigcup_{n<\omega} M_n$.

The same applies to limit stages in resolutions for our MAECs.

% \textcolor{red}{This is the reason because when we consider
%   resolutions in MAECs, we have to take the completion of the union of
%   an $\ordencea$-increasing chain of models in $\C{K}$ instead of just
%   its union as in discrete AECs.}

\begin{remark}\label{UnionCountable}
Let $\langle M_i : i<\sigma \rangle$ be an $\ordencea$-increasing chain
of models in an MAEC $\C{K}$. If $b\in
\overline{\bigcup_{i<\sigma}}M_i$, there exists a Cauchy sequence
$(a_n)_{n<\omega}$ in $\bigcup_{i<\sigma} M_i$ such that
$(a_n)_{n<\omega}\to b$. If $cf(\sigma)>\omega$, clearly $\overline{
  \bigcup_{i<\sigma} M_i}= \bigcup_{i<\sigma} M_i$. However, if
$cf(\sigma)=\omega$, we may well have that $b\in \overline{
  \bigcup_{i<\sigma} M_i}\setminus  \bigcup_{i<\sigma} M_i$.
% \\
% \indent If $cf(\sigma)>\omega$, there exists some $j<\sigma$ such that
% for every $k<\omega$ we have that $a_k\in M_j$, since $M_j$ is
% complete then $b\in M_j\subseteq \bigcup_{i<\sigma} M_i$ (because $b$
% is a limit of a sequence of elements of $M_j$), hence $\overline{
%   \bigcup_{i<\sigma} M_i}= \bigcup_{i<\sigma} M_i$.
% \\
% \indent If $cf(\sigma)=\omega$, given $\varepsilon>0$ there exists
% some $k<\omega$ such that $d(a_k,b)<\varepsilon$ (because
% $(a_n)_{n<\omega}\to b$). Since $a_k\in  \bigcup_{i<\sigma} M_i$, so
% $d(b, \bigcup_{i<\sigma} M_i)<\varepsilon$ (by
% definition~\ref{dist_element_set}, because $a_k\in \bigcup_{i<\sigma}
% M_i$) but not necessarily $b\in \bigcup_{i<\sigma} M_i$.
\end{remark}

\begin{definition}\label{k_embed}
We call a function $f:M\to N$  a $\mathcal{K}$-\emph{embedding} if
\begin{enumerate}
\item For every $k$-ary function symbol $F$ of $L$, we have $f(F^M(a_1\cdots
  a_k))=F^N(f(a_1)\cdots f(a_k))$.
\item For every constant symbol $c$ of $L$, $f(c^M)=c^N$.
\item For every $m$-ary relation symbol $R$ of $L$, for every
  $\bar{a}\in M^{m}$,
  $d(\bar{a},(R^M)^{-1}[k])=d(f(\bar{a}),(R^N)^{-1}[k])$ for all
  $k\in [0,1]$ (where $d(c,X):=\inf \{d(c,b): b\in X\}$, for  $c\in
  M$ and $X\subset M$, or for arities above $1$, using the product metric).
\item $f[M]\ordencea N$.
\end{enumerate}
\end{definition}

\begin{definition}[Amalgamation Property, AP]\label{AP}
%\new{10-03-12}
Let $\C{K}$ be an MAEC. We say that $\C{K}$ satisfies {\it Amalgamation Property} (for short {\it AP}) if and only if  for every $M,M_1,M_2\in \C{K}$, if $g_i:M\to M_i$ is a $\C{K}$-embedding $(i\in \{1,2\})$ then there exist $N\in \C{K}$ and $\C{K}$-embeddings $f_i:M_i\to N$ ($i\in \{1,2\}$) such that $f_1\circ g_1=f_2\circ g_2$.
\end{definition}

\[
\begin{diagram}
\node{M_1} \arrow{e,t,..}{f_1} \node{N}\\
\node{M} \arrow{n,l}{g_1} \arrow{e,b}{g_2} \node{M_2} \arrow{n,r,..}{f_2}
\end{diagram}
\]

\begin{definition}[Joint Embedding Property, JEP]
%\new{10-03-12}
Let $\C{K}$ be an MAEC. We say that $\C{K}$ satisfies {\it Joint Embedding Property} (for short {\it JEP}) if and only if  for every $M_1,M_2\in \C{K}$ there exist $N\in \C{K}$ and $\C{K}$-embeddings $f_i:M_i\to N$ ($i\in \{1,2\}$).
\end{definition}

\begin{remark}[Monster Model]\label{Monster_Model}
%\new{10-03-12}
If $\C{K}$ is an MAEC which satisfies AP and JEP and has large enough models, then we can
construct a large enough model $\monster$ (which we call a {\it
  Monster Model}) which is homogeneous --i.e., every isomorphism
between two $\C{K}$-substructures of $\monster$ can be extended to an
automorphism of $\monster$-- and also universal --i.e., every model
with density character $<dc(\monster)$ can be $\C{K}$-embedded into
$\monster$.
\end{remark}

\begin{definition}[Galois type]\label{Galois_Type}
%\new{10-03-12}
Under the existence of a monster model $\monster$ as in
remark~\ref{Monster_Model}, for all $\tuple{a}\in \monster$ and
$N\ordencea \monster$, we define $\gatp(\tuple{a}/N)$ (the
\emph{Galois type of $\tuple{a}$ over $N$}) as the orbit of
$\tuple{a}$ under $Aut(\monster/N):=\{f\in Aut(\monster): f\rest N =
id_{N}\}$. We denote the space of Galois types over a model $M\in
\C{K}$ by $\gaS(M)$.
\end{definition}

Throughout this paper, we assume the existence of a homogenous and
universal monster model as in remark~\ref{Monster_Model}.

\begin{definition}[Distance between types]\label{Distance_Types}
%\new{10-03-12}
Let $p,q\in \gaS(M)$. We define $d(p,q):=\inf\{d(\tuple{a},\tuple{b}):
\tuple{a},\tuple{b}\in \monster, \tuple{a}\models p , \tuple{b}\models
q\}$, where ${\rm lg}(\tuple{a})={\rm lg}(\tuple{b})=:n$ and
$d(\tuple{a},\tuple{b}):=\max\{d(a_i,b_i): 1\le i\le n\}$.
\end{definition}

It is worth noting here that this distance between types (defined
originally by Hirvonen in her thesis~\cite{Hi} but really generalizing
usual distances between orbits) gives a pseudo-metric, and makes the
set of types into a pseudo-metric space. It is easier for technical
purposes to actually get a \emph{metric on the set of types}. This is
the purpose of our next definition (Continuity of Types), a property
we will assume throughout the rest of the paper. Notice that in
principle we may have two Galois types at a small distance, but still
have many realizations between them being at possibly arbitrarily
large distances between them.

\begin{definition}[Continuity of Types]\label{Continuity_Types}
%\new{10-03-12}
Let $\C{K}$ be an MAEC and consider $(a_n)\to a$ in $\monster$. We say
that $\C{K}$ satisfies {\it Continuity of Types
  Property}\footnote{This property is also called {\it Perturbation
    Property in \cite{Hi}}} (for short, {\it CTP}), if and only if, if
$\gatp(a_n/M)=\gatp(a_0/M)$ for all $n<\omega$ then
$\gatp(a/M)=\gatp(a_0/M)$.
\end{definition}

\begin{remark}
%\new{10.03.12}
In general, distance between types $d$ (see
Definition~\ref{Distance_Types}) is just a pseudo-metric. But it is
straightforward to see that the fact that $d$ is a metric is
equivalent to CTP.
\end{remark}

Throughout this paper, we will assume our MAECs satisfy the CTP (so,
distance between types is in fact a metric).

\section{Limits and Independence Notions}
\label{sec:limitindep}

The basic tools of our proof - its statement, the control of the
construction - are given in this section. First of all, we provide a
warm-up, simple proof of the existence of limit models in
Galois-stable MAECs; along the proof, some definitions used later will
be provided. Later, we provide the reader a reminder of two notions of
independence which were defined and studied in detail in~\cite{ViZa1}:
$\varepsilon$-splitting ($\varepsilon$-independence) and smooth
independence. We will briefly describe their origin and state the main
results, referring the reader to the detailed proofs in the reference
provided. 

\subsection{Existence of limit models in Galois-stable MAECs}
We now adapt one of the existing notions of limit models 
(see \cite{GrVa}) to the metric context. 
This still leaves open the many variants of the notion that have
recently been used by Shelah in NIP contexts (see for
instance~\cite{Sh900}), or in strictly stable first order contexts.

\begin{definition}[Universality]\label{Universality}
%\new{10.03.12}
Let $\C{K}$ be an MAEC with CTP and $N\ordencea M$. We say that $M$ is
$\lambda$-d-{\it universal} over $N$ iff
for every $N'\ordenceac N$ with density character $\lambda$ there
exists a $\C{K}$-embedding $f:N'\to M$ such that $f\rest N=id_N$. We
say that $M$ is {\it universal} over $N$ if $M$ is $dc(M)$-universal
over $N$.
\end{definition}

The following lemma will be useful later: it provides relative
saturation criteria by iterating $\omega$-many times dense relative
saturation.

\begin{lemma}\label{saturation}
Suppose that we have an increasing $\ordencea$-chain of models
$(N_n:n<\omega)$ such that $N_{n+1}$ realizes a dense subset of
$\gaS(N_n)$. Then, every type in $\gaS(N_0)$ is realized in
$N_\omega:=\overline{\bigcup_{n<\omega} N_n}$.
\end{lemma}

\bdem
See \cite[Lemma 1.19]{ViZa1}.
\edem

%PZ junio 7: UNIVERSALIDAD YA SE HABÍA DEFINIDO ARRIBA.
%\begin{definition}
%Let $M,N\in \C{K}$ be such that $M\ordencea N$. We say that $N$ is
%$\mu$-$\mathbf{d}$-universal over $N$ iff for every $N'\ordenceac N$
%such that $dc(N)=\mu$ we have that there exists a $\ordencea$-embedding
%$f:N'\to M$ which fixes pointwise $M$. We say that $N$ is
%$\mathbf{d}$-universal over $M$ iff $N$ is
%$dc(M)$-$\mathbf{d}$-universal over $M$.
%\end{definition}

We drop the prefix $\mathbf{d}$ if it is clear that we are working in
a metric setting.

\begin{definition}[Limit models]\label{LimitModels_def}
Let $M,N\in \C{K}$ be such that $M\ordencea N$, where $dc(M)=\mu$.
We say that $N$ is $(\mu,\theta)$-$\mathbf{d}$-limit over $M$ iff
there exists an increasing and continuous $\ordencea$-chain
$(M_i:i<\theta)$ such that $M_0:=M$, $\overline{\bigcup_{i<\theta}M_i}=N$,
$dc(M_i)=\mu$ for every $i<\theta$ and also $M_{i+1}$ is
$\mu$-$\mathbf{d}$-universal over $M_i$.
\end{definition}

\begin{definition}
We say that $\C{K}$ is $\mu$-$\mathbf{d}$-stable iff for every
$M\in\C{K}$ such that $dc(M)\le \mu$ we have that
$dc(\gaS(M))\le\mu$
\end{definition}

We now prove the existence of universal extensions in the setting of
Metric Abstract Elementary Classes. We point out that this is an
adaptation of the proof of the existence of universal extensions
over a given model $M$ in the setting of Abstract Elementary Classes
(see \cite{GrVa}). In that proof, under $\mu$-stability, we can
consider an increasing and continuous $\C{K}$-chain $\langle
M_i:i<\mu\rangle$ such that $M_0:=M$ and where $M_{i+1}$ realizes
every Galois-type in $\gaS(M_i)$. So, $\bigcup_{i<\mu}M_i$ is
universal over $M$. But in this setting, we cannot consider directly
from $\mu$-$\mathbf{d}$-stability that $M_{i+1}$ realizes every type
in $\gaS(M_i)$. But we use~Lemma~\ref{saturation} in a suitable way
for guaranteeing that requirement.

\begin{proposition}[Existence of universal extensions]\label{Existence_Universal}
Let ${K}$ be an MAEC $\mu$-$\mathbf{d}$-stable with AP and CTP. Then
for
all ${M}\in \C{K}$ such that $dc({M})=\mu$ there exists
${M}^{*}\in \C{K}$ universal over ${M}$. such that
$dc({M}^*)=\mu$
\end{proposition}

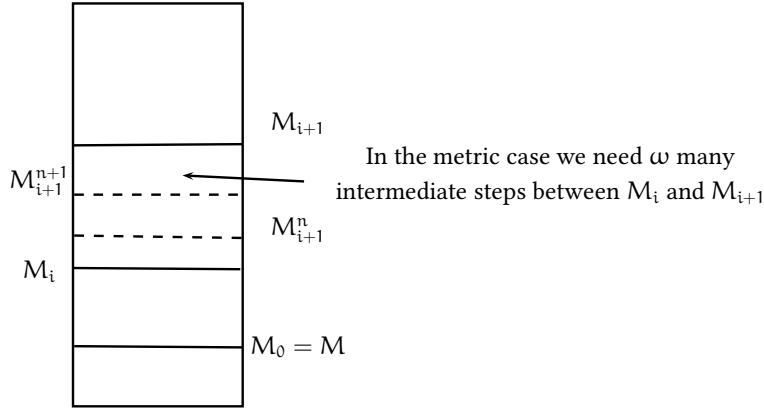
\scalebox{.8} % Change this value to rescale the drawing.
{
\begin{pspicture}(0,-3.37)(19.062813,3.37)
\psframe[linewidth=0.04,dimen=outer](3.9009376,3.37)(1.0409375,-3.37)
\psline[linewidth=0.04cm](1.0809375,-2.35)(3.8809376,-2.37)
\psline[linewidth=0.04cm](1.0609375,-1.05)(3.8409376,-1.07)
\psline[linewidth=0.04cm](1.0609375,0.99)(3.8809376,1.01)
\psline[linewidth=0.04cm,linestyle=dashed,dash=0.16cm 0.16cm](1.0809375,0.17)(3.8009374,0.17)
\psline[linewidth=0.04cm,linestyle=dashed,dash=0.16cm 0.16cm](1.0809375,-0.51)(3.8609376,-0.55)
\psline[linewidth=0.04cm,arrowsize=0.05291667cm 2.0,arrowlength=1.4,arrowinset=0.4]{->}(4.9009376,0.39)(2.8809376,0.49)
\rput(9,0.78){In the metric case we need $\omega$ many}
\rput(9,0.20){intermediate steps between $M_i$ and $M_{i+1}$}
\rput(4.8,-2.36){$M_0=M$}
\rput(0.5,-1.1){$M_i$}
\rput(4.8,-0.42){$M_{i+1}^n$}
\rput(0.5,0.4){$M_{i+1}^{n+1}$}
\rput(4.8,1.34){$M_{i+1}$}
\end{pspicture}
}

\bdem
The proof follows almost along the same lines as
the proof of existence
of universal models in usual AECs (see~Claim~2.9~of~\cite{GrVa} and
Claim~1.15.1 of~\cite{Sh600}); that is, by trying to capture realizations of types
along the construction in a coherent way, and building the universal
extension as a union of a chain (we do not repeat all the details of
the proof, but point out the differences).

In our metric setting, we need to be careful with the way we
realize the types along the construction: although this cannot be done in an
immediate way in each successor stage as in~\cite{GrVa},
lemma~\ref{saturation} provides the realizations we need of dense
subsets of the typespace in $\omega$ many steps.

We construct an increasing and continuous $\ordencea$-chain of
models $({M}_i:i<\mu \rangle$ such that ${M}_0:={M}$,
${M}_{i+1}$ is the completion of the union of a resolution
$({M}_{i+1}^{n}:n<\omega)$ where ${M}_{i+1}^0:={M}_i$, 
${M}_{i+1}^{n+1}$ realizes a dense subset of $\gaS({M}_{i+1}^n)$
and $dc({M}_{i+1}^n)=\mu$ for every $n<\omega$. This is possible by
$\mu$-$\mathbf{d}$-stability of $\C{K}$. Take
${M}^*:=\overline{\bigcup_{i<\mu} {M}_i}$.  ${M}^*$ turns out to be
universal over ${M}$ --- by the same argument as
in~Claim~2.9~of~\cite{GrVa}.
\edem[Prop. \ref{Existence_Universal}]

\begin{corollary}\label{Existence_Limit_Models}
Let $\C{K}$ be an MAEC $\mu$-$\mathbf{d}$-stable with AP. Then for
all ${M}\in \C{K}$ such that $dc({M})=\mu$ there exists
${M}^{*}\in \C{K}$ limit over ${M}$ such that
$dc({M}^*)=\mu$
\end{corollary}
\bdem Iterate the construction given in proposition
\ref{Existence_Universal}.
\edem[Cor. \ref{Existence_Limit_Models}]

\subsection{Smooth independence in MAECs.}

We define the notion of {\it smooth
  independence} for MAECs with CTP and AP, and state some of its
properties. For a more complete
analysis of this independence notion, see~\cite{ViZa1,Za}. In
particular, our paper~\cite{ViZa1} contains proofs of the properties of
both smooth independence and $\varepsilon$-splitting. We quote the
results that appeared in~\cite{ViZa1} here since we need them for our
proof in what remains of the paper.

The idea of $\varepsilon$-splitting and $\eindep$ is a variant of
splitting for abstract elementary classes - a notion widely and
productively used in the vast literature on the subject since the
early work of Shelah. For many purposes, until now, splitting (well,
rather, non-splitting) has been the most robust notion of independence
in AECs - unlike the situation in first order logic, where in many
ways non-forking has played that role. Here we adapt the definition
first to a fixed $\varepsilon$ - saying that two types in the
definition are ``at least $\varepsilon$-apart from one another''
instead of just saying they are different. This gives
$\varepsilon$-splitting and $\eindep$. Then we take a bolder step
(necessary for dealing with chains of uncountable cofinality later):
we define ``smooth independence'' as a way of ``blending'' all
$\varepsilon$-independences, \emph{along a resolution}. Smooth
independence is the
most robust property for our purposes, but we will use both
$\varepsilon$-non-splitting independence and smooth independence
in our Assumption~\ref{Superstability}. The reader is urged to refer
to~\cite{ViZa1} for the proofs of the results quoted here.

\begin{definition}[$\varepsilon$-splitting and $\eindep$]
Let $N\ordencea M$ and $\E>0$. We say that $\gatp(a/M)$
$\E$-\emph{splits} over $N$ iff there exist $N_1,N_2$ with
$N\ordencea N_1,N_2\ordencea
M$ and $h:N_1\approx_N N_2$ such that
$\dist(\gatp(a/N_2),h(\gatp(a/N_1))\ge \E$. We use  $a\eindep_N
M$ to denote the fact that  $\gatp(a/M)$ does
not $\varepsilon$-split over $N$.
\end{definition}

\begin{definition}\label{r_independence}
Let $N\ordencea M$. Fix $\mathcal{N}:=(N_i : i<\sigma \rangle$
a resolution of $N$. We say that $a$ is \emph{smoothly independent}
from $M$ over $N$ relative to $\C{N}$ (denoted by $a\indep^{\C{N}}_N
M$) iff for every $\E>0$ there exists $i_\E<\sigma$
such that $a\eindep_{N_{i_\E}} M$.
\end{definition}

% We call {\it smooth independence} the notion of independence given
% above, inspired by \cite{BaSh}. In that paper, J. Baldwin and
% S. Shelah defined {\it smoothness} as a nice property of an abstract
% class of models $\C{K}$ which involves increasing chains of models,
% context where the existence of a kind of monster model holds.

% In the following lines, we provide a list of properties -without
% proofs- of smooth independence. For the proofs, see \cite{ViZa1,Za}.

\begin{proposition}[Monotonicity of
  non-$\varepsilon$-splitting]\label{monotonicity_eps}
Let $M_0\ordencea M_1\ordencea M_2\ordencea M_3$. If $a\eindep_{M_0}
M_3$ then $a\eindep_{M_1} M_2$.
\end{proposition}

\begin{proposition}[Monotonicity of
  smooth-independence]\label{monotonicity}
Let $M_0\ordencea M_1\ordencea M_2\ordencea M_3$. Fix
$\mathcal{M}_k:=(M^k_i : i<\sigma_k \rangle$ a resolution of
$M_k$ ($k=0,1$), where $\C{M}_0\subseteq \C{M}_1$. If
$a\indep^{\C{M}_0}_{M_0} M_3$ then $a\indep^{\C{M}_1}_{M_1} M_2$.
\end{proposition}

\begin{lemma}[Stationarity (1)]\label{Tarea1}
Suppose that $N_0\ordencea N_1 \ordencea N_2$ and $N_1$ is
universal over $N_0$. If $\gatp(a/N_1)=\gatp(b/N_1)$,
$a\eindep_{N_0} N_2$ and $b\eindep_{N_0} N_2$, then
$\dist(\gatp(a/N_2),\gatp(b/N_2)<2\varepsilon$.
\end{lemma}

%\begin{proposition}[Extension of $\indep^{\C{N}}$ over universal
%models]\label{ExtUniv}\index{Materias}{Extension Property}
%If $N\ordencea M \ordencea M'$, $\C{N}:=(N_i:i<\sigma \rangle$
%is a resolution of $N$, $M$ is universal over $N$ and
%$p:=\gatp(a/M)\in\gaS(M)$ is a Galois type such that
%$a\indep^{\C{N}}_N
%M$, then there exists $b$ such that $\gatp(b/M)=\gatp(a/M)$ and
%$b\indep^{\C{N}}_N M'$.
%\end{proposition}

\begin{proposition}[Stationarity (2)]\label{UnExtUniv}
If $N\ordencea M \ordencea M'$, $M$ is universal over $N$,
$\C{N}:=(N_i:i<\sigma \rangle$ a resolution of $N$ and
$p:=\gatp(a/M)\in\gaS(M)$ is a Galois type such that
$a\indep^{\C{N}}_N
M$, then there exists an unique extension $p^*\supset p$ over $M'$
which is independent (relative to $\C{N}$) from $M'$ over $N$.
\end{proposition}

\begin{proposition}[Continuity of
  smooth-independence]\label{cont_indep}
Let $(b_n)_{n<\omega}$ be a convergent sequence and
$b:=\lim_{n<\omega} b_n$. If $b_n\indep^{\C{N}}_N M$ for every
$n<\omega$, then $b\indep^{\C{N}}_N M$.
\end{proposition}

\begin{proposition}[stationarity (3)]\label{Stationarity}
Let $M_0\ordencea M\ordencea N$ be such that $M$ is a
$(\mu,\sigma)$-limit model over $M_0$, witnessed by
$\C{M}:=(M_i:i<\sigma\rangle$. If $a,b\indep^{\C{M}}_M N$ and
$\gatp(a/M)=\gatp(b/M)$, then we have $\gatp(a/N)=\gatp(b/N)$.
\end{proposition}

\begin{proposition}[Transitivity]
Let $M_0\ordencea M_1\ordencea M_2$ be such that $M_0$ is a limit
model over some $M'\ordencea M_0\ordencea M_1$ (witnessed by
$\C{M}_0$) and $M_1$ is a limit model over $M_0$ (witnessed by
$\C{M}_1'$). Let $\C{M}_1:=\C{M}_0	\ ^\frown \C{M}_1'$ (i.e., $\C{M}_1$ is the concatenation, as sequences, of $\C{M}_0$ and $\C{M}_1'$), %so $\C{M}_0\subset \C{M}_1$.
Then $a\indep^{\C{M}_0}_{M_0} M_2$ iff
$a\indep^{\C{M}_0}_{M_0} M_1$ and $a\indep^{\C{M}_1}_{M_1} M_2$.
\end{proposition}

The three versions of stationarity correspond to three situations
(over universal models for $\varepsilon$-non-splitting and smooth
independence, and over limit models) and provide the technical tools
to guarantee (in the first case) \emph{small distance} between
non-splitting
extensions. Of course, in first order or even in discrete AECs, the
word ``stationarity'' refers to \emph{uniqueness} of non-forking (or
non-splitting) extensions under various assumptions (stability). Here,
the assumptions are given over the base model (universality or being a
limit - implicitly, Galois stability) and the consequences are not
quite uniqueness in the first case, but \emph{small distance}. We get
that
$\varepsilon$-non-splitting or smooth independent extensions of a type
are \emph{very close} in the metric of types. (Notice however that this
distance between orbits does not imply that arbitrary realizations of
the types are close.) The proofs of these forms of stationary are in
Lemma 3.6, Proposition 3.9 and Proposition 3.15 of~\cite{ViZa1}.
\\ \\
\indent The following property of smooth independence (called {\it
  anti-reflexivity}) is a metric version of the following property of
forking in the first order setting: if $a\indep^{\th}_B a$ then
$a\in acl(B)$. We use this basic property in proving that
the completion of the union of a $\le$-increasing chain of towers is
again a tower (fact~\ref{DefCompl}).

\begin{proposition}[Anti-reflexivity]\label{s_antireflexivity}
Let $M\ordencea N$ where $M$ is a $(\mu,\theta)$-limit model witnessed
by $\C{M}:=\{M_i:i<\theta\}$. If $a\indep_{M}^{\C{M}}N$ and $a\in N$,
then $a\in M$.
\end{proposition}

As mentioned above, the proofs of the propositions in this section are in
our paper~\cite{ViZa1}.

\section{Towers: the core of the proof}
\label{sec:towers}

In this long section we provide the core of the proof - the end of
proof is presented in the next, final, section. Here, the tools of the
construction (smooth towers, reduced towers, the full relativeness of
smooth towers) are introduced, and the main lemmas of the construction
are proved. For the sake of making the proof more readable, we have
split this section into three subsections, dealing with smooth towers,
$d$-reduced towers and full relativeness.

Throughout this section, we assume that all our models have density
character $\mu$, all orderings denoted by $I,I',I_\beta$, etc. have
cardinality $\mu$ as well.
%AÒadido PZ, marzo 15 de 2010
%and \textcolor{red}{$cf(I)=cf(I')=cf(I_\beta)>\omega$, unless we state the
%contrary}\marginpar{HAY QUE QUITAR ESTA FRASE - PROBLEM\'ATICA. PZ.}.
Furthermore, we only consider models which are not compact
as metric spaces (in first order continuous model theory, compact
metric spaces are the analogue of finite models in first order
discrete model theory - by ruling out compact metric spaces, we rule
out in advance that exceptional case).

\begin{assumption}[superstability]\label{Superstability}
For every $a$ and every increasing and continuous $\ordencea$-chain of
models
$(M_i : i<\sigma \rangle$ and $\C{M}_j$ a resolution of $M_j$
($j<\sigma$):
\begin{enumerate}
\item If $p\rest M_i \indep^{\C{M}_0}_{M_0} M_i$ for all $i<\sigma$, then
$p\indep^{\C{M}_0}_{M_0} \overline{\bigcup_{i<\sigma} M_i}$.

\item if $cf(\sigma)>\omega$, there exists $j<\sigma$ such that
$a\indep^{\C{M}_j}_{M_j} \bigcup_{i<\sigma} M_i$.
\item if $cf(\sigma)=\omega$, for every $\varepsilon>0$ there exists
  $j<\sigma$ such that
$a\indep^{\E}_{M_j} \overline{\bigcup_{i<\sigma} M_i}$.
\end{enumerate}
\end{assumption}

We need to break into two cases, according to cofinality, our
treatment of independence above: under countable cofinality, an
element $b$ belonging to the completion of a union is witnessed by a
Cauchy sequence $a_i\to b$ from the union, and we get the
approximation to our independence directly from these elements
$a_i$. Under uncountable cofinality, $b$ has to belong to the union,
and we transfer the approximation idea to resolutions.

\begin{remark}\label{TwoImpliesThree}
\end{remark}
Notice that if $cf(\sigma)>\omega$, (2) in our assumption above gives
$j<\sigma$ such that $a\indep^{\C{M}_j}_{M_j} \bigcup_{i<\sigma}
M_i=\overline{\bigcup_{i<\sigma} M_i}$ (equality by
Remark~\ref{UnionCountable}). Smooth independence then gives, for each
positive integer $n$, an index $j_n$ such that
$a\indep^{\frac{1}{n}}_{M_{j_n}} \bigcup_{i<\sigma} M_i$ (with
$M_{j_n}$ in the resolution $\C{M}_j$ of $M_j$). Monotonicity
of $\indep^{\frac{1}{n}}$ (proposition~\ref{monotonicity_eps}) then
yields $a\indep^{\frac{1}{n}}_{M_{j}} \bigcup_{i<\sigma} M_i$. So, the
conclusion of (3) in our Assumption~\ref{Superstability} holds also if
$cf(\sigma)>\omega$.

\begin{remark}\label{whysuperstability}
\end{remark}
Why do we use the blanket term ``superstability''? Let us notice that
in the first order case, our hypothesis holds for superstable
theories.
Let $T$ be a superstable (discrete) first order theory. Suppose an
$\prec$-increasing and continuous chain of models $\langle M_i:
i<\omega_1 \rangle$ and let $b$ any element in a monster model of
$T$. By superstability, there exists a finite $A\subset
\bigcup_{i<\omega_1} M_i$ such that $tp(b/\bigcup_{i<\omega_1} M_i)$
does not fork over $A$. Since there exists some $j<\omega_1$ such that
$A\subset M_j$, by monotonicity of non-forking we also have that
$tp(b/\bigcup_{i<\omega_1} M_i)$ does not fork over $M_j$. In stable
theories, non-forking agrees with non-splitting. Since we are in a
discrete setting, considering $\varepsilon\le1$
non-$\varepsilon$-splitting agrees with non-splitting.

\begin{remark}\label{rem-categoricity}
\end{remark}
It is now worth explaining a further crucial point concerning the use
of categoricity in our proof. Historically, the most desirable result
for uniqueness of limit models would have only used the
``superstability assumptions'', as in first order logic. However, the
situation in Abstract Elementary Classes in general seems to hold
additional complexity: not only does one have to use categoricity to
obtain the symmetry of non-splitting (clearly not immediate from
definition, and at the moment not even clearly a \emph{consequence} of
uniqueness of limit models), but there is the additional fact that
categoricity assumptions do not seem to imply our superstability
assumptions in the metric AEC case. We do have of course that for
reasonable enough AECs and MAECs, categoricity at a cardinal $\lambda$
implies Galois-stability above $LS(\C{K})$ and below $\lambda$ (as long
as there is enough in the AEC to get Ehrenfeucht-Mostowski models, a
proof due to Hyttinen and explained in~\cite[8.20 and 8.21]{BaMon} provides
stability). Superstability-like properties seem harder to obtain in
our metric AEC context.

Additionally, in the case of \emph{discrete} Abstract Elementary
Classes (see \cite[Ch. 15]{BaMon}) extensive use of the non-metric
version of our Assumption~\ref{Superstability} has been derived from
categoricity, and used. When our metric is discrete, our MAECs are
AECs, and our Assumption corresponds to superstability assumptions of
that setting. Our notions are adapted to the metric context in ways
that both reflect the original AEC setting and respond to the more
sophisticated continuity   % The additional complication arising in
                           % our setting
% corresponds to the need to separate chains of cofinality $\omega$ from
% chains of uncountable cofinality in our proofs.

\subsection{Smooth towers.}

\begin{definition}[smooth towers]%\index{Materias}{Smooth tower}
Let $I$ be a well-ordering, $\mathfrak{M}:=(M_i:i\in I)$ be an
$\ordencea$-increasing chain, $\overline{a}:=(a_i:i\in I)$,
$\mathfrak{N}:=(N_i:i<\sigma)$ be a sequence of models in $\C{K}$,
$\C{M} :=({\C{M}_j}: j\in I )$ be a sequence of resolutions ${\C{M}
  _j}$ of
$M_j$ ($j\in I$) and $\C{N} :=({\C{N}_j}: j\in I)$ be a sequence of
resolutions ${\C{N}_j}$ of $N_j$ ($j\in I$). We
say that $\tower$ is an \emph{$I$-smooth tower} (shortly, {\it $I$-s-tower})
iff for every $i\in I$ we have that
$M_i$ is a $(\mu,\sigma)$-limit model over $N_i$ for some
$\sigma<\mu^+$, $a_i\in
M_{i+1}\setminus M_i$ and $a_i\indep^{\C{N}_i}_{N_i} M_i$. If the set of indexes is clear, we may drop it. Also, if it is clear that we are working in a metric context, by a {\it tower} we mean an $s$-tower.%We denote the class of $I$-smooth-towers where all itsby $\C{K}_I$}.
\end{definition}

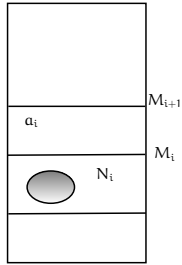
\begin{center}
\scalebox{.5} % Change this value to rescale the drawing.
{
\begin{pspicture}(0,-3.47)(5.2690625,3.47)
\psframe[linewidth=0.04,dimen=outer](3.72,3.47)(0.0,-3.47)
\psline[linewidth=0.04cm](0.04,-2.15)(3.7,-2.13)
\psline[linewidth=0.04cm](0.02,-0.59)(3.68,-0.57)
\psline[linewidth=0.04cm](0.04,0.71)(3.68,0.71)
\rput(4.2,-0.56){$M_i$}
\rput(4.2,0.86){$M_{i+1}$}
\rput(0.64921874,0.28){$a_i$}
\psellipse
[linewidth=0.04,dimen=outer,fillstyle=gradient,gradlines=2000,gradbegin=gray,gradend=white,gradmidpoint=1.0](1.17,-1.44)(0.65,0.45)
\rput(2.6079688,-1.1){$N_i$}
\end{pspicture}
}
\end{center}

Roughly speaking, an s-tower is composed of a $\ordencea$-increasing
(not necessarily continuous) chain of models $\M{M}:=(M_i:i\in I)$ and
a collection of $\C{K}$-submodels $\M{N}:=(N_i:i\in I)$ such that each
$M_i$ is a $(\mu,\sigma)$-limit model over $N_i$ (for some
$\sigma<\mu^+$) which codify a smooth independence of the elements
$a_i$ taken in the s-tower (i.e., $a_i\indep_{N_i}^{{\C{N}_i}}M_i$).

\begin{definition}[Extension of
  s-towers]\label{ExtTowers}
Let $I\le I'$ be well-orderings, %\linebreak
%\marginpar{\textcolor{red}{PZ:jun8}}
$\tower$ %\in \C{K}_{\mu,I}$ PZ: junio 8 2014
be an $I$-s-tower and 
and $(\mathfrak{M}',\overline{a}',\mathfrak{N}',M ', N ')$%\in \C{K}_{\mu,I'}$. PZ junio 8 2014.
be an $I'$-s-tower.
We say that %\linebreak
$(\mathfrak{M}',\overline{a}',\mathfrak{N}',M ', N ')$ extends
$\tower$  (denoted by %\linebreak
$(\mathfrak{M}',\overline{a}',\mathfrak{N}',M ', N ')\ge \tower$)
iff  for every $i\in I$:%%\new{}  
\begin{enumerate}
    \item $M_i'$ is a proper universal model over $M_i$
    \item ${\C{M} _i}$ is an initial segment of ${\C{M} _i'}$, as sequences.
    \item $a_i=a_i'$
    \item $N_i=N_i'$
    \item ${\C{N}_i}={\C{N}_i'}$
\end{enumerate}
\end{definition}

%ExtensiÛn de torres - eliminado julio 8 de 2011
\begin{center}
% Generated with LaTeXDraw 2.0.8
% Wed Aug 17 20:07:22 COT 2011
% \usepackage[usenames,dvipsnames]{pstricks}
% \usepackage{epsfig}
% \usepackage{pst-grad} % For gradients
% \usepackage{pst-plot} % For axes
\scalebox{0.7} % Change this value to rescale the drawing.
{
\begin{pspicture}(0,-3.34)(9.2,3.36)
\definecolor{color96b}{rgb}{0.6,0.6,0.6}
\psframe[linewidth=0.04,dimen=outer](2.44,2.36)(0.0,-3.24)
\psframe[linewidth=0.04,dimen=outer](7.28,2.38)(3.36,-3.28)
\psline[linewidth=0.04cm](5.34,2.32)(5.38,-3.32)
\psline[linewidth=0.04cm](0.02,-2.14)(2.46,-2.14)
\psline[linewidth=0.04cm](3.36,-2.12)(7.3,-2.08)
\psline[linewidth=0.04cm](0.08,-0.98)(2.38,-1.0)
\psline[linewidth=0.04cm](3.4,-0.98)(7.22,-0.96)
\psline[linewidth=0.04cm](0.0,0.18)(2.38,0.18)
\psline[linewidth=0.04cm](3.36,0.22)(7.28,0.2)
\usefont{T1}{ptm}{m}{n}
\rput(1.8,-1.6){$M_i$}
\usefont{T1}{ptm}{m}{n}
\rput(2.0,-0.155){$M_{i+1}$}
\psellipse[linewidth=0.04,linestyle=dashed,dash=0.16cm 0.16cm,dimen=outer,fillstyle=solid,fillcolor=color96b](0.84,-1.62)(0.6,0.42)
\psellipse[linewidth=0.04,linestyle=dashed,dash=0.16cm 0.16cm,dimen=outer,fillstyle=solid,fillcolor=color96b](4.12,-1.65)(0.6,0.37)
\usefont{T1}{ptm}{m}{n}
\rput(4.2,-1.6){$N_i$}
\usefont{T1}{ptm}{m}{n}
\rput(0.8,-1.595){$N_i$}
\psline[linewidth=0.04cm,fillcolor=color96b,linestyle=dashed,dash=0.16cm 0.16cm,doubleline=true,doublesep=0.12](5.38,-0.7)(7.22,-0.7)
\psline[linewidth=0.04cm,fillcolor=color96b,linestyle=dashed,dash=0.16cm 0.16cm,doubleline=true,doublesep=0.12](5.36,-0.34)(7.2,-0.34)
\psline[linewidth=0.04cm,fillcolor=color96b,linestyle=dashed,dash=0.16cm 0.16cm,doubleline=true,doublesep=0.12](5.36,0.02)(7.24,0.0)
\psline[linewidth=0.04cm,fillcolor=color96b,linestyle=dashed,dash=0.16cm 0.16cm,doubleline=true,doublesep=0.12](5.36,-1.78)(7.26,-1.8)
\psline[linewidth=0.04cm,fillcolor=color96b,linestyle=dashed,dash=0.16cm 0.16cm,doubleline=true,doublesep=0.12](5.36,-1.44)(7.18,-1.46)
\psline[linewidth=0.04cm,fillcolor=color96b,linestyle=dashed,dash=0.16cm 0.16cm,doubleline=true,doublesep=0.12](5.36,-1.12)(7.22,-1.14)
\psdots[dotsize=0.12](3.62,-0.14)
\psdots[dotsize=0.12](0.26,-0.14)
\usefont{T1}{ptm}{m}{n}
\rput(0.64,-0.275){$a_i$}
\usefont{T1}{ptm}{m}{n}
\rput(4.12,-0.235){$a_i$}
\usefont{T1}{ptm}{m}{n}
\rput(4.5,-0.75){$M_{i+1}$}
\usefont{T1}{ptm}{m}{n}
\rput(8.0,-0.175){$M_{i+1}'$}
\usefont{T1}{ptm}{m}{n}
\rput(8.0,-1.395){$M_i'$}
%\psline[linewidth=0.04cm,fillcolor=color96b,arrowsize=0.05291667cm 2.0,arrowlength=1.4,arrowinset=0.4]{->}(1.6,2.74)(4.56,2.7)
%\usefont{T1}{ptm}{m}{n}
%\rput(3.08,3.165){$f$}
\end{pspicture} 
}
\end{center}

\subsection{$d$-reduced towers.}

%Nueva definicion propuesta por PZ - enero 6 de 2015
\begin{definition}\label{reduced_towers}
Let $\delta:\mathbb{R}^+\rightarrow \mathbb{R}^+$
be a mapping such that for every $\alpha,\E>0$
\begin{enumerate}
\item $\delta(\alpha\cdot \E) = \alpha\cdot \delta(\E)$
  and
\item $\delta(\varepsilon)\le \varepsilon$
\end{enumerate}

%The last condition ($\delta(\varepsilon)\le \varepsilon$)
Notice that for instance the function
$\delta(\varepsilon)=\varepsilon$ satisfies the two properties
above. The function $\delta$ will provide, given $\varepsilon$ a
``continuity modulus'' for the definition of reduced towers adapted to
our context.

An $I$-s-tower $\tower$ is said to be a {\it $d$-reduced tower
  (relative to $\delta$)} iff
whenever $\tower '\ge \tower$ then for
every $j\in I$ and every $\E>0$, if $b\in \overline{\bigcup_{i\in I}
  M_i}$ and $d(b,M'_j)<\delta(\E)$ then $d(b,M_j)<\E$.
\end{definition}

\begin{notation}\label{delta_epsilon}
Let $\delta:\mathbb{R}^+\rightarrow \mathbb{R}^+$ be a mapping as in
definition~\ref{reduced_towers}. For the sake of simplicity, we denote
$\delta_\varepsilon:=\delta(\varepsilon)$.
\end{notation}

%Comentado por PZ enero 6 de 2015
%\begin{definition}
%An $I$-s-tower $\tower$ is said to be a {\it $d$-reduced tower} iff there exists a mapping $\delta:\mathbb{R}^+\rightarrow \mathbb{R}^+$ such that $\delta(\alpha\cdot \E)=\alpha\cdot \delta(\E)$ for every $\alpha,\E>0$, in such a way that if $\tower '\ge \tower$ then for every $j\in I$ and every $\E>0$, if $b\in \overline{\bigcup_{i\in I} M_i}$ and $d(b,M'_j)<\delta_\E:=\delta(\E)$ then $d(b,M_j)<\E$.
%\end{definition}

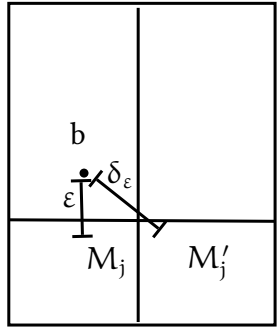
\begin{center}
% Generated with LaTeXDraw 2.0.8
% Fri Mar 23 15:06:00 COT 2012
% \usepackage[usenames,dvipsnames]{pstricks}
% \usepackage{epsfig}
% \usepackage{pst-grad} % For gradients
% \usepackage{pst-plot} % For axes
\scalebox{1} % Change this value to rescale the drawing.
{
\begin{pspicture}(0,-2.65)(3.58,2.69)
\psframe[linewidth=0.04,dimen=outer](3.58,1.67)(0.0,-2.65)
\psline[linewidth=0.04cm](1.74,1.59)(1.74,-2.59)
\psline[linewidth=0.04cm](0.02,-1.23)(3.56,-1.25)
\rput(1.32,-1.8){$M_j$}
\rput(2.67,-1.8){$M'_j$}
\psdots[dotsize=0.12](1.02,-0.61)
\psline[linewidth=0.04cm,tbarsize=0.07055555cm 5.0]{|-|}(1.16,-0.67)(2.04,-1.37)
\psline[linewidth=0.04cm,tbarsize=0.07055555cm 5.0]{|-|}(0.98,-0.69)(1.0,-1.47)
\rput(1.5,-0.625){$\delta_\E$}
\rput(0.82,-0.965){$\E$}
\rput(0.94,-0.085){b}
\rput(0.3,2.2){$\tower$}
\rput(3.2,2.2){$\tower'$}
\end{pspicture} 
}

\end{center}
%\marginpar{PZ:jun8. Añadí esa línea porque el referee dijo que el dibujo no era claro.}
In this picture, the left column corresponds to the tower $\tower$ and the right one corresponds to the tower $\tower'$.

%\marginpar{\textcolor{red}{PZ:jun8. NUEVO, a petición del referee.}}
\begin{definition}[Completion of the union of an increasing chain of
  towers]%Directed limit of an increasing chain of towers.
Let $\langle I_\alpha: \alpha<\beta \rangle$ be an $\le$-increasing
chain of well-orderings\footnote{Here, the $I_\alpha$ are indeed
  \emph{arbitrary} well-orderings; for this definition we do not
  demand the well-orderings to have uncountable cofinality.} and
suppose that $I:=\bigcup_{\alpha<\beta}
I_\alpha$ is a well-ordering. Let $\langle \tower^\alpha :
\alpha<\beta \rangle$ be an $\le$-increasing chain of towers, where
$\tower^\alpha$ is an $I_{\alpha}$-s-tower. %We define the
                                %\emph{completion of the union} of
                                %this $\le$-increasing chain as
                                %follows: %\new{} 
\begin{enumerate}
	\item For every $i\in I$, let $\alpha(i):=\min\{\alpha<\beta:
          i\in I_\alpha\}$.
	\item For $i\in I$, %define $M:=\lim_{\rightarrow} \langle
                            %\overline{\bigcup_{i\in
                            %I_\alpha}M^\alpha_i} ;f_{\alpha,\gamma}:
                            %\alpha\le \gamma<\beta \rangle$, with
                            %canonical embeddings denoted by
                            %$f_{\alpha,\beta}$. 
	define $M_i:=\overline{\bigcup_{\alpha(i)\le \alpha<\beta}
          M^\alpha_i}$, and $N_i:=N^{\alpha(i)}_i$ and
        $a_i:=a^{\alpha(i)}_i$ as they were defined in the
        $I_{\alpha(i)}$-s-tower $\tower^{\alpha(i)}$. Define
        $\mathfrak{M}:=(M_i:i\in I)$, $\overline{a}:=(a_i:i\in I)$,
        $\mathfrak{N}:=(N_i:i\in I)$, $\C{M}:=({\C{M}_i}:i\in I)$
        where each ${\C{M}_i}$ is defined as the respective
        concatenations of $\C{M}^{\alpha(i)}_i$ and
        ${\C{M}^{\alpha+1}_{i}}\setminus \C{M}^{\alpha}_{i}$ for
        $\alpha(i)\le \alpha<\beta$ (remember that since
        $\tower^{\alpha'}\le \tower^{\alpha}$ for $\alpha(i)\le
        \alpha'<\alpha<\beta$, then $\C{M}^{\alpha'}_{i}$ is an
        initial segment of $\C{M}^{\alpha}_{i}$, by definition of
        $\le$) and $\C{N}:=({\C{N}^{\alpha(i)}_i}:i\in I)$.

The $I$-s-tower $\tower$ defined as above is called the {\it
  completion of the union} of the sequence\linebreak
$\langle \tower^\alpha :\alpha<\beta\rangle$.
\end{enumerate} 

% Generated with LaTeXDraw 2.0.8
% Sun Jun 08 16:17:42 COT 2014
% \usepackage[usenames,dvipsnames]{pstricks}
% \usepackage{epsfig}
% \usepackage{pst-grad} % For gradients
% \usepackage{pst-plot} % For axes
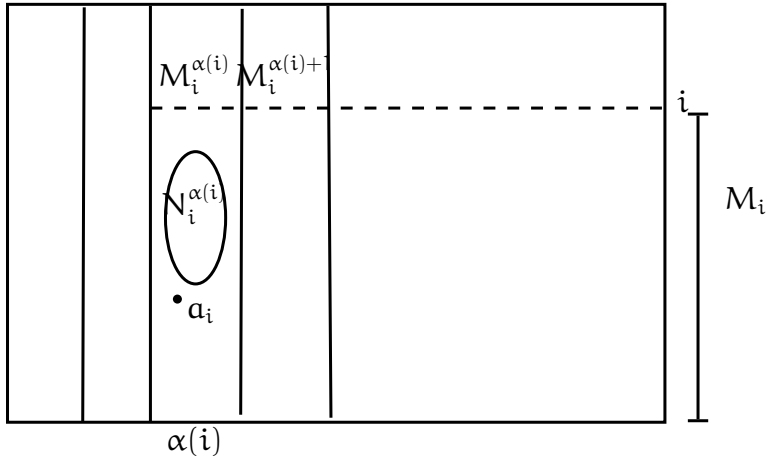
\begin{center}
\scalebox{1} % Change this value to rescale the drawing.
{
\begin{pspicture}(0,-2.988496)(12.121895,2.948496)
\psframe[linewidth=0.04,dimen=outer](8.78,2.948496)(0.0,-2.6515038)
\psline[linewidth=0.04cm](1.04,2.8884962)(1.02,-2.6315038)
\psline[linewidth=0.04cm](1.92,2.9084961)(1.92,-2.6315038)
\psline[linewidth=0.04cm](3.14,2.8684962)(3.12,-2.531504)
\psline[linewidth=0.04cm](4.28,2.9084961)(4.32,-2.5715039)
%\psline[linewidth=0.04cm,linestyle=dashed,dash=0.16cm 0.16cm](1.92,0.3884961)(8.72,0.34849608)
\usefont{T1}{ptm}{m}{n}
\rput(2.5,-2.9){$\alpha(i)$}
\psline[linewidth=0.04cm,linestyle=dashed,dash=0.16cm 0.16cm](1.92,1.5484961)(8.72,1.5484961)
\usefont{T1}{ptm}{m}{n}
\rput(2.5,2){$M^{\alpha(i)}_i$}
\usefont{T1}{ptm}{m}{n}
\rput(3.7,2){$M^{\alpha(i)+1}_i$}
\psellipse[linewidth=0.04,dimen=outer](2.52,0.0884961)(0.42,0.9)
\usefont{T1}{ptm}{m}{n}
\rput(2.5,0.2734961){$N^{\alpha(i)}_i$}
\psdots[dotsize=0.12](2.28,-0.9915039)
\usefont{T1}{ptm}{m}{n}
\rput(2.581455,-1.1865039){$a_i$}
%\psellipse[linewidth=0.04,dimen=outer](2.45,0.8984961)(0.39,0.35)
%\usefont{T1}{ptm}{m}{n}
%\rput(2.4,0.9){$N_i$}
%\psdots[dotsize=0.12](2.84,1.3084961)
%\usefont{T1}{ptm}{m}{n}
%\rput(3.4,1.2){$a_i$}
%\usefont{T1}{ptm}{m}{n}
\rput(9,1.633496){$i$}
\psline[linewidth=0.04cm,tbarsize=0.07055555cm 5.0]{|*-|}(9.2,1.4684961)(9.2,-2.6315038)
\usefont{T1}{ptm}{m}{n}
\rput(9.832627,0.31349608){$M_i$}
\end{pspicture} 
}
\end{center}
In this picture, each column corresponds to a tower. Therefore, the column in the right side of each tower corresponds to an $\le$-extension as s-towers.
\\

%\marginpar{\textcolor{red}{PZ:jun8. NUEVO, a petición del referee.}}
\begin{fact}\label{DefCompl}
$\tower$ is in fact an s-tower.
\end{fact}
\bdem
Notice that $M_i$ is universal over $N_i$ (since $M^{\alpha(i)}_i$ is
universal over $N_i$ and $M_i^\alpha(i)\ordencea M_i$),
$a_i\indep^{{\C{N}_i}}_{N_i} M^\alpha_i$ for every $\alpha(i)\le
\alpha<\beta$ (by definition of s-tower), so by superstability
(assumption \ref{Superstability})  $a_i\indep^{{\C{N}_i}}_{N_i} M_i$
(since $M_i:=\overline{\bigcup_{\alpha(i)\le \alpha<\beta}
  M^\alpha_i}$).

It is worth mentioning here that we use anti-reflexivity
(Proposition~\ref{s_antireflexivity}) implicitely in these
constructions: If $a_i\indep^{{\C{N}_i}}_{N_i} M^\alpha_i$ and
$a_i\notin N_i$ then $a_i\notin M^\alpha_i$.

%\textcolor{red}{
By definition of extension of s-towers (definition~\ref{ExtTowers}
(1.)) $M_i^{\gamma+1}$ is universal over $M_i^{\gamma}$ ($\alpha(i)\le
\gamma<\beta$, $\beta$ a limit ordinal). So,
$M_i:=\overline{\bigcup_{\alpha(i)\le \alpha<\beta} M^{\alpha}_i}$ is
a limit model over $M^{\alpha(i)}_i$. Since by definition of s-tower
$M^{\alpha(i)}_i$ is a limit model over $N_i=N^{\alpha(i)}_i$,
concatenating the sequences witnessing $M_i$ is a limit model over
$M^{\alpha(i)}_i$ and $M^{\alpha(i)}_i$ is a limit model over $N_i$
respectively, we have that $M_i$ is a limit model over
$N_i$. %}%\marginpar{REVISAR- PZ JULIO 1 2014}

The tower $\tower$ defined as above is in fact a tower which
$\le$-extends every tower $\tower^\alpha$: Let $i\in I$, since
$M^{\alpha(i)+1}_i$ is universal over $M^{\alpha(i)}_i$ (since
$\tower^{\alpha(i)}\le \tower^{\alpha(i)+1}$, by definition of $\le$)
and $M_i\ordenceac M^{\alpha(i)+1}_i$, then $M_i$ is universal over
$M^{\alpha(i)}_i$.

\edem[Fact \ref{DefCompl}]
\end{definition}

\begin{proposition}[Density of reduced towers]\label{reduced_density}
Fix $\delta:\mathbb{R}^+\to \mathbb{R}^+$ a mapping as in
definition~\ref{reduced_towers}.
Given  $\tower$ an $I$-s-tower, there exists an $I$-s-tower $\tower'
\ge \tower$ such that $\tower'$ is a $d$-reduced tower relative to
$\delta$.
\end{proposition}
\bdem
Suppose the proposition is false. Let $(\tower^{\alpha}:
\alpha<\mu^+ \rangle$ be an $\le$-increasing chain of towers such that
$\tower^{\alpha+1}\ge \tower^{\alpha}$ witnesses that
$\tower^{\alpha}$ is not a reduced tower relative to
  $\delta$,
where $\tower^0:=\tower$.
\\ \\
Let $\tower^{\mu^+}$ be the completion of the union of %\newline
 $\langle
\tower^{\alpha}: \alpha<\mu^+ \rangle$. Fix
$\E>0$ and $b\in \overline{\bigcup_{i\in I} M_i^{\mu^+}}$; define 
$i_\E(b):=\min\{ i\in I: \text{\ there exists } b'\in \bigcup_{j\le i}
M^{\mu^+}_i$ such that  $d(b,b')<\E \}$ and 
$\zeta_\E(b):=\min \{\zeta<\mu^+: \text{\ there exists } b'\in
M^{\zeta}_{i_\E(b)} \text{ such that } d(b',b)<\E \}$. Let
$E:=\{\zeta<\mu^+:$ for all $b\in \overline{\bigcup_{i\in I}
  M_i^{\zeta}},$ $\zeta_{\E}(b)<\zeta \}$. $E$ is a stationary subset
of $\mu^+$, in particular non-empty. Let $\zeta\in E$. 

By our construction, $\tower^{\zeta+1}$ witnesses that $\tower^\zeta$ is not a
d-reduced tower relative to $\delta$.
There exist therefore some $\E>0$  
%for every linear mapping $\delta:\mathbb{R}^+\to \mathbb{R}^+$ there exists
and $b\in \overline{\bigcup_{i\in I} M^\zeta_i}$ such that for some level
$j\in I$ we have that 
\begin{eqnarray}\label{EQN}
d(b,M^{\zeta+1}_j)<\delta(\E) \text{\ but \ }%$d(b,M^\zeta_j)= PZ junio 8 2014
d(b,M^\zeta_j)\ge \E .
\end{eqnarray}
%but $d(b,M^\zeta_j)= PZ junio 8 2014 $d(b,M^\zeta_j)\ge \E$.
%Since $b\in \overline{\bigcup_{i\in I} M^\zeta_i}$, there exists a sequence  $(b_n)$ in $\bigcup_{i\in I} M^\zeta_i$ such that $(b_n)\to b$. So, there exists $N<\omega$ such that $d(b_N,b)<\E/2$. %Let $j\in I$ such that $b_N\in M_j^{\zeta}$, therefore $i_{\E/2}(b)\le j$. Notice that $\zeta_{\E/2}(b)<\zeta$, since $\zeta\in E$. Since $b_N\in M_{i_{\E/2}(b)}^{\zeta}$, there exists $(c_n)\in {\bigcup_{\alpha<\zeta}M^\alpha_{i_{\E}(b)}}$ such that $(c_n)\to b_N$. Let $M<\omega$ be such that $d(c_M,b_N)<\E/2$, so $d(c_M,b)\le d(c_M,b_N)+d(b_N,b)<\E$. Let $\alpha<\zeta$ be such that $c_M\in M^\alpha_{i_{\E/2}(b)}$. Notice that $d(b,M^\zeta_{j})\le d(b,M^\zeta_{i_{\E/2}(b)})\le d(b,M^\alpha_{i_{\E/2}(b)})<\E$ 	(contradicts (1) above).
Since $d(b; M^{\mu^+}_{j})\le
d(b;M^{\zeta+1}_j)<\delta(\varepsilon)\le \varepsilon$ (by property
2. of $\delta$ in definition~\ref{reduced_towers}), then
$i_{\varepsilon}(b)\le j$.
Since $\zeta\in E$, we have that $\zeta_{\varepsilon}(b)<\zeta$;
therefore $d(b;M^{\zeta}_j)\le d(b;M^{\zeta}_{i_{\varepsilon}(b)})\le
d(b;M^{\zeta_{\varepsilon}(b)}_{i_{\varepsilon}(b)})<\E$ But this contradicts
(\ref{EQN}) above.
\edem[Prop. \ref{reduced_density}]

\begin{proposition}\label{Union_reduced}
Let $\delta:\mathbb{R}^+\to \mathbb{R}^+$ be a function as in
definition~\ref{reduced_towers} and $(\tower_i : i<\beta \rangle$ be
an $\le$-increasing
sequence of $d$-reduced towers relative to $\delta$. Then the
completion of its union is a
$d$-reduced tower relative to some $\delta'$ as in
definition~\ref{reduced_towers}.
\end{proposition}
\bdem
Let $\tower_\beta$ be the completion of the union of $\langle
\tower_i
: i<\beta \rangle$ and let \linebreak
$\tower'\ge \tower_\beta$.
\\ \\
Let $\E>0$ and $K<\omega$ be such that ${2\delta_{\E}}/{K}<\E$
(remember that by notation~\ref{delta_epsilon},
  $\delta_\E:=\delta(\varepsilon)$). 
\\ \\
Let $\E':=\E-2\delta_{\E}/K>0$ and
$\E'':=\min\{\delta_{5\E'/6}/10,\delta_{\E}/K\}$. 
\\ \\
Let $b\in \overline{\bigcup_{i\in I} M_i^\beta}$ be such that
$d(b;M_j')<\delta_{5\E'/6}$. Since $b\in \overline{\bigcup_{i\in I}
  M_i^\beta}$, then there exists a sequence $(a_n)\in \bigcup_{i\in I}
M_i^\beta$ such that $(a_n)\to b$. So, there exists $n_\E<\omega$ such
that 
\begin{eqnarray*}
d(a_{n_\E},b)<\E''
\end{eqnarray*}
Let $k_\E\in I$ be such that $a_{n_\E}\in M_{k_\E}^\beta$. By
definition of $\tower_\beta$,
$M_{k_{\E}}^\beta:=\overline{\bigcup_{\alpha(k_\E)\le \alpha<\beta}
  M_{k_\E}^\alpha}$, hence there exists $(c_{m})\in
\bigcup_{\alpha(k_\E)\le \alpha<\beta}M_{k_\E}^\alpha$ such that
$(c_m)\to a_{n_\E}$. Therefore, there exists ${m_\E}<\omega$ such that 
\begin{eqnarray*}
d(c_{m_\E},a_{n_\E})<\E''
\end{eqnarray*}
Let $\alpha_\E<\beta$ be such that $c_{m_\E}\in
M_{k_\E}^{\alpha_\E}$. %Let $d_{m_\E}\in  M_{k_\E}^{\alpha_\E}$ be
                       %such that
                       %$c_{m_\E}=f_{\alpha_\E,\beta}(d_{m_\E})$. %Suppose
                       %that $d(b;M_j^\beta)\ge \E$, hence $d(b,x)\ge
                       %\E$ for every $x\in M_j^\beta$.
\\ \\
Let $y\in M_j'$ be such that $d(b,y)<\delta_{5\E'/6}$ (this exists
because $d(b,M_j')<\delta_{5\E'/6}$). Notice that 
\begin{eqnarray*}
d(c_{m_\E},y) &\le& d(c_{m_\E},a_{n_\E})+d(a_{n_\E},b)+d(b,y)\\
&<&  2\E'' + \delta_{5\E'/6}\\
&\le& 2\delta_{5\E'/6}/10 + \delta_{5\E'/6}\\
&=& \frac65\delta_{5\E'/6}\\
&=& \delta_{\frac65 \cdot \frac{5\E'}{6}}\\
&=& \delta_{\E'}
\end{eqnarray*}
Since  $\tower'\ge \tower_{\alpha_\E}$ and 
\begin{eqnarray*}
d(c_{m_\E},y)&<&\delta_{\E'}\\
&=&{\delta_{\E-2\delta_{\E}/K}}
\end{eqnarray*}
then $d(c_{m_\E};M_j')<\delta_{\E-2\delta_{\E}/K}$; since
$\tower_{\alpha_\E}$ is $d$-reduced relative to $\delta$, then
$d(c_{m_\E};M_{j}^{\alpha_\E})\le {\E-2\delta_{\E}/K}$.
\\ \\
Let $x\in M_j^{\alpha_\E}$ be such that
$d(c_{m_\E},x)<\E-2\delta_{\E}/K$. Notice that
\begin{eqnarray*}
d(b,x) &\le& d(b,a_{n_\E})+d(a_{n_\E},c_{m_\E})+ d(c_{m_\E},x)\\
&<& 2\E''+(\E-2\delta_{\E}/K)\\
&\le& 2\delta_{\E}/K+ (\E-2\delta_{\E}/K)\\
&=& \E
\end{eqnarray*}
Since $M_j^{\alpha_\E}\ordencea M_j^\beta$,  then
$d(b,M_j^\beta)<\E$. Notice that $\tower_\beta$ is a $d$-reduced tower
relative to $\delta':\mathbb{R}^+\to \mathbb{R}^+$ defined as
$\delta'(\varepsilon):=\delta(5\E/6)$.
\edem[Prop. \ref{Union_reduced}]

%\textcolor{red}{
%\begin{remark}\label{delta_prime}
%Let $\delta:\mathbb{R}^+\to \mathbb{R}^+$ be a mapping as in definition~\ref{reduced_towers} and define $\delta':\mathbb{R}^+\to \mathbb{R}^+$ as in the end of the proof of proposition~\ref{Union_reduced}. Notice that for every $\E>0$
%\begin{eqnarray*}
%\delta'(\E) &=& \delta (5\E/6)\\
%&=& 5/6\cdot \delta(\E)\\
%&<& 1\cdot \delta(\E)\\
%&=& \delta(\E)
%\end{eqnarray*}
%Therefore, $\delta'(\E)<\delta(\E)$ for every $\E>0$.
%\end{remark}
%}

\begin{fact}\label{SatIso}
If $M,N\in \C{K}$ are $\lambda$-d-saturated models of density
character $\lambda^+$ then $M\approx N$.
\end{fact}
\bdem
By a standard back and forth argument, as in discrete AECs.
\edem[Fact \ref{SatIso}.]

The following proof is central in the argument (continuity of
d-reduced towers of the correct density character) and is the place
where we use our categoricity (in $\mu^+$) assumption. Our use of the
categoricity assumption has been adapted from VanDieren's proof
in~\cite{Va2}.

\begin{proposition}\label{reduced_continuity}
Let $\C{K}$ be $\mu^+$-categorical, and let $\delta$ be a
function as in~\ref{reduced_towers}.  Then every d-reduced
(relative to $\delta$) tower of density character $\mu$
is continuous.
\end{proposition}
\bdem
The plan of the proof is to assume failure of continuity
  of a d-reduced tower, toward the construction of a model of density
  character $\mu^+$ that omits some type - this will contradict the
  categoricity of the class in $\mu^+$ since, as we will describe, the
  class also must have a Galois-saturated model of that density
  character. The type omitted in the model constructed will be
  obtained from our assumption of failure of continuity.

Suppose then that there is a d-reduced tower $\tower$ which is not
continuous; i.e.,
there exists a limit element $\delta\in I$ and $b\in M_\delta$ such
that $b\notin \overline{\bigcup_{i<\delta} M_i}$. By the density
property (Prop. \ref{reduced_density}), combined with the
non-maximality of our models, we can build a strictly
$\le$-increasing sequence of d-reduced towers $(\tower_\alpha:
\alpha<\mu^+ \rangle$ such that $\tower_0:=\tower\rest \delta$. Let 
$\tower_{\mu^+}:= \overline{\bigcup_{\alpha<\mu^+}
  \tower_{\alpha}}$.\\
Our non-Galois saturated model of density character
  $\mu^+$ is $M:=\overline{\bigcup_{i<\delta}
  M_i^{\mu^+}}$. Indeed, let us show that $M$ omits the Galois type
$p=\gatp(b/\overline{\bigcup_{i<\delta} M_i})$. Assume toward a
contradiction that $M$ realizes $\gatp(b/\overline{\bigcup_{i<\delta}
M_i})$. So, there are an element $c\in M$ and a map $f\in
Aut(\monster/\overline{\bigcup_{i<\delta} M_i})$ such that
$f(c)=b$. Define 
$\tower ' :=f[\tower_{\mu^+}]$.\\
 Notice that
$\tower_0=f[\tower_0]\le \tower '$ and $b=f(c)\in \overline{\bigcup_{i<\delta}
  M_i'}$. %(remember that by hypothesis, $b\in \overline{\bigcup_{i<\delta}
  %M_i'}$).%\marginpar{\textcolor{red}{PZ:jun8, petición del referee.}}
 Let $\E:=d(b;\overline{\bigcup_{i<\delta}M_i})> 0$ (as we have
 $b\notin \overline{\bigcup_{i<\delta}M_i}$). Since
$b\in \overline{\bigcup_{i<\delta} M_i'}$ there exists $(b_n)\in
\bigcup_{i<\delta} M_i'$ such that $(b_n)\to b$. We can assume that
there exists $N<\omega$ such that $d(b_N,b)<\delta(\E)/2$. Let
$i<\delta$ be such that $b_N\in M_i'$. Since
$M_i'=f[\overline{\bigcup_{\alpha<\mu^+}
  M_i^\alpha}]=\overline{\bigcup_{\alpha<\mu^+} f[M_i^\alpha]}$,
there exists a sequence $(c_n)\in {\bigcup_{\alpha<\mu^+}
  f[M_i^\alpha]} $ such
that $(c_n)\to b_N$. We can find $K<\omega$ such that
$d(c_K,b_N)<\delta(\E)/2$ and some $\alpha<\mu^+$ such that $c_K\in
f[M^{\alpha}_i]$,%\marginpar{\textcolor{red}{PZ:jun8}}
 so $d(c_K,b)\le
d(c_K,b_N)+d(b_N,b)<\delta(\E)/2+\delta(\E)/2=\delta(\E)$. Notice that
$f[\tower_\alpha]\ge
\tower_0$ and $d(b;f[M_i^\alpha])\le
d(b,c_K)<\delta(\E)$. By the definition of reduced towers,
  since
$\tower$ is a d-reduced tower we have
$d(b;\overline{\bigcup_{i<\delta}M_i})\le
d(b;M_i)<\E$. But this contradicts our initial choice of
  $\E$ as $d(b;\overline{\bigcup_{i<\delta}M_i})$.

Therefore, $M$ is a model in the class of density character
$\mu^+$ which is not $\mu$-d-saturated. By using Ehrenfeucht-Mostowski
models, there also exists a model
$N\in \C{K}$ $\mu$-d-saturated, of density character $\mu^+$. Therefore
$M\not\approx N$, contradicting our $\mu^+$-categoricity assumption.
\edem[Prop. \ref{reduced_continuity}]

  The previous proof exhibits an interesting analogy to existing
  proofs of stability below categoricity cardinals: a contradiction to
  categoricity is obtained by assuming failure of stability, and
  building two models, one omitting some types, the second one
  saturated. In our previous proof, categoricity in $\mu^+$ entails
  \emph{continuity of d-reduced towers in density character $\mu$,}
  not directly stability. One may thus regard continuity of these towers as a
  (very localized) form of ``stability''.

\subsection{Full-relativeness of s-towers}

\begin{definition}[strong type]\label{StrongType}
Let $M$ be a $\sigma$-limit model
\begin{enumerate}
\item $\stype(M):=\left\{
(p,N) :
\begin{tabular}{l}
$N\ordencea M$\\
$N$ is a $(\mu,\theta)$-limit model\\
$M$ is universal over $N$\\
$p\in \gaS(M)$ is non-algebraic\\
and $p\indep^{\C{N}}_{N} M $\\
for some resolution $\C{N}$ of $N$.
\end{tabular}\right\}$

\item Two strong types $(p_l,N_l)\in \stype(M_l)$ ($l\in
\{1,2\}$) are {\it parallel} (which we denote by
$(p_1,N_1)\parallel (p_2,N_2$) iff for every $M'\ordenceac
{M}_1,{M}_2$ with density character $\mu$, there exists $q\in
\gaS(M')$ which extends both $p_1$ and $p_2$ and
$q\indep^{\C{N}_l}_{N_l} M'$ ($l\in\{1,2\}$)
(where $\C{N}_l$ is the resolution of $N_l$ which satisfies
$p_i\indep^{\C{N}_l}_{N_l} M_l$).
\end{enumerate}
\end{definition}

\begin{assumption}\label{In}
  \begin{itemize}
  \item Through this subsection, assume that $I_0$ is a well ordered
    set which
    % has a cofinal sequence $(i_\alpha:\alpha<\theta)$
    has a sequence $(i_\alpha:\alpha\le \theta)$ where ${\rm
      cf}(\theta)>\omega$, $i_\theta$ is the maximum of $I_0$ and
    $(i_\alpha:\alpha<\theta)$
    is a cofinal sequence in $I\setminus \{i_\theta\}$. Moreover we will
    use the following
    notation. In the proofs of Weak Full Relativeness
    (\ref{weak_full}) below and in the proof
    of~\ref{omega_theta} wrapping up the large scale proof, we use
    orderings $I_n$ for $n<\omega$. %$I_0$ is some fixed well-ordering
                                %ofcofinality $\theta$ and QUIT… esto
                                %porque ya se aclarÛ que es $I_0$
                                %atr·s.
  \item $I_{n+1}$ is constructed from $I_n$ by
  inserting a new $\mu$-sequence between each pair of elements
  $i_\alpha$ and its $I_0$-successor $i_{\alpha +1}$ ($\alpha<\theta$)
  in such a way
  that in the ordering $I_n$ there are $n$ many ``copies'' of
  $\mu$ between $i_\alpha$ and $i_{\alpha+1}$. This notation
  corresponds to the analog notation in a proof for discrete AECs
  in~\cite{GrVaVi}.
  \end{itemize}
\end{assumption}

\begin{definition}[Metric s-Towers]
An s-tower $\tower$ is called a \emph{metric s-tower} if the
resolution witnessing that  $M_i$ is a $(\mu,\sigma)$-limit model over
$N_i$ is spread-out. A spread-out resolution $ \C{M}$ of $M$ is a
resolution where for every $\gamma$, $M^{\gamma+1}$ is an
$\omega_1$-limit over $M^\gamma$.
%% la idea es que las resoluciones esparcidas suben a lo largo de
%% lÌmites, y las s-torres mÈtricas suben  a los lrgo de resoluciones
%% esparcidas
\end{definition}

\begin{definition}[Full relative s-towers]
Let $\tower$ be a s-tower indexed by $I$. Let
$(M_i^\gamma:\gamma<\sigma)$ be a sequence which witnesses that
$M_i$ is a $(\mu,\sigma)$-limit model. We say that $\tower$ is a
relative s-tower with respect to $(M_i^{\gamma})_{i\in
I,\gamma<\sigma}$ iff for every $i_\alpha\le i<i_{\alpha+1}$ and
$(p,M_i^\gamma)\in \stype(M_i)$ there exists $i\le j<i_{\alpha+1}$
such that $(p,M_i^\gamma)\parallel (\gatp(a_j/M_j),N_j)$.
\end{definition}

\begin{proposition}\label{limit1}
Suppose that for every $\alpha<\theta$ there are $\mu\cdot \omega$
many elements between $i_\alpha$ and $i_{\alpha+1}$. Let $\tower$ be a
full relative s-tower with respect to $(M_i^{\gamma})_{i\in
  I,\gamma<\sigma}$. Then $M:=\overline{\bigcup_{i\in I} M_i}$ is a
limit model over $M_{i_0}$.
\end{proposition}
\bdem It is enough to prove that $M_{i_{\alpha+1}}$ is universal over $M_
{i_\alpha}$.  Let $p:=\gatp(a/M_{i_\alpha})\in
\gaS(M_{i_\alpha})$ and $\E>0$. So, by assumption \ref{Superstability} (3.)
there exists $\gamma:=\gamma_{\E}<\sigma$ such that
$a\eindep_{M^{\gamma_\E}_{i_0}} M_{i_0}$ (if $cf(\sigma)=\omega$
we have to apply part 3. of assumption~\ref{Superstability},  if
$cf(\sigma)>\omega$, this follows from assumption~\ref{Superstability}
(2.) and remark~\ref{TwoImpliesThree}).
\\ \\
%AÒadido marzo 7 de 2010
By construction, $M^{\gamma+1}_{i_\alpha}$ is a $(\mu,\omega_1)$-limit
model over $M^{\gamma}_{i_\alpha}$.
%
%Eliminado marzo 7 de 2010
%$M^{\gamma}_{i_0}\ordencea N\ordencea M^{\gamma+1}_{i_0}$ such
%that $N$ is an $\omega_1$-limit model over $M^{\gamma}_{i_0}$.
Let $(M^*_i:i<\omega_1)$ be a resolution which witnesses that.
\\ \\
Consider $q:=p\upharpoonright M^{\gamma+1}_{i_\alpha}$, so by assumption
\ref{Superstability} (2.) there exists $i<\omega_1$ such that
$q\indep^{\C{M}^*_i}_{M^*_i} M^{\gamma+1}_{i_\alpha}$
(as the cofinality is uncountable, we apply
  part 2. of assumption~\ref{Superstability}). By extension over
universal
models~(proposition~\ref{UnExtUniv},  %\cite[2.7]{ViZa1}
notice that $M^{\gamma+1}_{i_\alpha}$
is universal over $M^*_i$), there exists $q^*\in\gaS(M_{i_\alpha})$ an extension
of $q$ such that $q^*\indep^{\C{M}^*_i}_{M^*_i} M_{i_\alpha}$. So, $(q^*,M^*_i)\in
\stype(M_{i_\alpha})$. By full relativeness of $\tower$, there exists
$i_\alpha\le j_1 <i_{\alpha+1}$ such that $(q^*,M^*_i)\parallel
(\gatp(a_j/M_j),N_j)$. Therefore, $q^*=\gatp(a_j/M_{i_\alpha})$ and so
$q^*$ is realized in $M_{j_1}$.
\\ \\
By monotonicity of non-$\E$-splitting, we have that $p$ does not
$\E$-split over $M^*_i$ (since $p$ does not $\E$-split over
$M^{\gamma}_{i_\gamma}$ and $M^{\gamma}_{i_\alpha}\ordencea M^*_i$); i.e.
$p\eindep_{M^*_i} M_{i_\alpha}$. Since $q^*\indep^{\C{M}^*_i}_{M^*_i}M_{i_\alpha}$, then
$q^*\eindep_{M^*_i} M_{i_\alpha}$ (by monotonicity of non-$\E$-splitting, proposition~\ref{monotonicity_eps}).
\\ \\
Also, since $q=p\upharpoonright M^{\gamma+1}_{i_\alpha}$ and $q^*\supset q$, then
$q^*\upharpoonright M^{\gamma+1}_{i_\alpha}= p\upharpoonright
M^{\gamma+1}_{i_\alpha}$. Notice that $M^{\gamma+1}_{i_\alpha}$ is
universal over $M^*_i$.
\\ \\
Since $p,q^*\indep^{\E}_{M^*_i} M_{i_\alpha} $, by a
weak version of stationarity (proposition~\ref{Tarea1}) we have that
$\dist(p,q^*)<2\E$.
Therefore, $M_{j_1}$ realizes a dense subset of $\gaS(M_{i_\alpha})$.
\\ \\
Doing a similar argument, we can construct an increasing sequence
$(j_n: n<\omega)$ in $I$ (where $j_0:=i_{\alpha}$) such that
$i_{\alpha}\le j_n <i_{\alpha+1}$, where $M_{j_{n+1}}$ realizes a
dense subset of $\gaS(M_{j_n})$.
\\ \\
Therefore, by lemma~\ref{saturation} we have that
$M^*:=\overline{\bigcup_{n<\omega} M_{j_n}}
\ordencea M_{i_{\alpha+1}}$ realizes every type over
$M_{j_0}=M_{i_\alpha}$, so $M_{i_{\alpha+1}}$
does.
\edem[Prop. \ref{limit1}]

The following fact is proved in a similar way as in the discrete case
(see \cite{GrVaVi}). For the sake of completeness, we give a proof of
this result.

\begin{proposition}\label{LimitModels}
If $\tower$
% \in \C{K}_{\mu,I}^* $
is an $I$-tower, there exists 
$(\mathfrak{M}',a,\mathfrak{N},\C{M}',\C{N})>\tower$
%in $\C{K}_{\mu,I}^*$
an $I$-tower such that for every limit $i\in I$, $M_i'$ is a
$(\mu,\mu)$-limit over
$\overline{\bigcup_{j<i}M_j}$
\end{proposition}

\bdem
First, we construct by induction on $i\in I$ a model $M_i^+\ordenceac
M_i$ and a directed system  $\langle f_{i,j}: i<j\in I\rangle $ of
$\ordencea$-embeddings (as in the discrete AEC case, one may prove
that the ``union axioms'' for metric AEC also hold for directed
systems) such that
$f_{i,j}:M_i^+\to M_j^+$ and $f_{i,j}\rest M_i=id_{M_i}$.
\\ \\
\indent Suppose $\langle M_k^+: k\le i\rangle$ and $\langle f_{k,l}: k<l\le i\rangle$ are
constructed. We give the construction of $M_{i+1}^+$ and
$f_{i,i+1}$. The construction of $f_{j,i+1}$ ($j<i$) are given by
definition of directed system. Let $M_{i+1}^*$ be a limit model over
$M_i^+$ and $M_{i+1}$. Since $a_{i+1}\indep^{\C{N}_{i+1}}_{N_{i+1}}
M_{i+1}$ and $M_{i+1}$ is universal over $N_{i+1}$ (by definition of
s-tower), by the extension property (proposition~\ref{UnExtUniv}) and
invariance of smooth independence there exists $f\in
Aut(\monster/M_{i+1})$ such that
$a_{i+a}\indep^{\C{N}_{i+1}}_{N_{i+1}}  f[M^*_{i+1}]$. Define
$M_{i+1}^+:=f[M_{i+1}^*]$ and $f_{i,i+1}:=f\rest M_{i}^+$.
\\ \\
\indent For limit $i\in I$, first take the directed limit of $\langle M_k^+:
k< i\rangle$ and $\langle f_{k,l}: k<l< i\rangle $  (that directed limit might
 not be contained in $\mathbb{M}$, a monster model in $\C{K}$,  but since $\mathbb{M}$ is universal and contains $M_i$ we may find
  $M_i^*\ordencea \mathbb{M}$  isomorphic to the directed limit of $\langle M_k^+:
k< i\rangle$ and $\langle f_{k,l}: k<l< i\rangle$ fixing $M_i$ pointwise)  and
then consider $M_{i}^+$ a limit model over this directed limit and
$(\mu,\mu)$-limit over $\overline{\bigcup_{j<i}f_{j,i}[M_j^+] }$.
\\ \\
\indent Fix $j\in I$. Let $f_{j,\sup(I)}$ and $M'_{j,\sup(I)}$ be
the respective directed limit of this directed system. Notice that
$f_{j,\sup(I)}\rest M_j=id_{M_j}$. Define
$M'_j:=f_{j,\sup(I)}[M_j^+]$.
\\ \\
\indent Notice that the s-tower
$(\mathfrak{M}',\tuple{a},\mathfrak{N},\C{M},\C{N})$ defined in this
way satisfies the requirements of the proposition.
\edem[Prop. \ref{LimitModels}]

\begin{lemma}[Weak Full Relativeness]\label{weak_full}
Given $\tower\in \C{K}_{\mu,I_n}^*$, there exists
$(\mathfrak{M}',a,\mathfrak{N},\C{M}',\C{N})>\tower$ in
$\C{K}_{\mu,I_{n+1}}^*$ such that for every $(p,N)\in \stype(M_i)$
(where $i\in I_n$ and $i_\alpha\le i<i_{\alpha+1}$)
there exists $i\le j < i_{\alpha+1}$ such that
$(\gatp(a_j/M_j'),N_j)\parallel (p,N)$.
\end{lemma}
\bdem
Let $M_{i_{\alpha+1}}'$ be a $(\mu,\mu)$-limit model over
$\overline{\bigcup_{j<i_{\alpha+1},j\in I_n} M_j}$ (by proposition
\ref{LimitModels}).
Let $(M_i' : l\in I_{n+i},i_\alpha+\mu\cdot n <l<{\alpha+1}
\rangle$ be an enumeration of
a resolution which witnesses that $M_{i_{\alpha+1}}'$ is
$(\mu,\mu)$-limit over the model $\overline{\bigcup_{j<i_{\alpha+1},j\in I_n}
  M_j}$.
\\ \\
\indent Let $\mathfrak{S}:=\{(p,N)^l_\alpha : i_\alpha +\mu\cdot n < l
< i_{\alpha+1}\}$ be an enumeration of a dense subset of
$\bigcup\{\stype(M_i): i\in I_n,i_\alpha\le
  i<i_{\alpha+1} \}$ (by $\mu$-stability).
Therefore, given $(p,N)^l_\alpha\in \mathfrak{S}$ there exists $i\in
I_n$ such that $i_\alpha\le i<i_{\alpha+1}$
such that $(p,N)^l_{\alpha}\in \stype{(M_i)}$. So $p^l_\alpha
\indep^{\C{N}^l_\alpha}_{N^l_\alpha} M_i$. Since
by definition of strong type $M_i$ is universal over $N^l_\alpha$ and
$M_i\ordencea M'_l$,
by~proposition~\ref{UnExtUniv} there exists
$p^*\in \gaS(M'_l)$ which extends $p^l_\alpha$ and
$p^*\indep^{\C{N}^l_\alpha}_{N^l_\alpha} M'_l$.
Notice that $M'_{succ_{I_{n+1}}(l)}$ is universal over $M'_l$ (by
construction), then there exists
$a_l\in  M'_{succ_{I_{n+1}}(l)}$ such that $a_l\models p^*$ (where $succ_{I_{n+1}}(l):=\min\{y\in I_{n+1}: l<y\}$; i.e., $succ_{I_{n+1}}(l)$ is the successor of $l$ in $I_{n+1}$). Consider
$N_l:=N^l_\alpha$.
So, $a_l\indep^{\C{N}_l}_{N_l} M'_l$. The s-tower constructed in this
way satisfies the requirements of the
proposition. \\ 
\edem[Lemma \ref{weak_full}]

\section{Uniqueness of Limit Models}

We now put together the material from the previous three sections and
finish the proof of uniqueness of limit models in the categorical case
for Metric Abstract Elementary Classes that have Amalgamation and
Continuity of Types (MAEC + AP + CTP).

Part of the outline of the proof is inspired in the proof of the
analogous results given by Grossberg, VanDieren and the first author
of this paper
in~\cite{GrVaVi}. There are however serious changes\emph{ in the
lemmas}, due to the difference in independence notions, in the
revised definition of ``reduced tower'' and in the proof of
continuity of reduced towers here. The metric context forces us to
thread finely and deal with differences that are not visible in the
usual (discrete) AEC context.

However, the results here follow a general outline of proof that
already has a history in the proof of Uniqueness of Limit Models in
``superstable'' AECs - in this very particular sense, this paper is a
contribution to the superstability of \emph{metric} Abstract
Elementary Classes where the types are orbital (AP) and are also
endowed with a metric (CTP).

% The following fact is inspired by the related result given in
% \cite{GrVaVi}. Although the sketch of the proof in the metric case is
% the same as the proof given in \cite{GrVaVi}, we have to point out
% that the details of the steps in the proof are quite
% different.%%\new{}   

\begin{proposition}\label{omega_theta}
Let $\C{K}$ be an MAEC satisfying AP, JEP, CTP, existence of large
enough models, assumption~\ref{Superstability} (superstability) and
which is $\mu^+$-categorical. There is a $(\mu,\theta)$-$d$-limit
model over $M$ which is also a $(\mu,\omega)$-$d$-limit model over
$M$.
\end{proposition}%%\new{}  
\bdem
If $cf(\theta)=\omega$, a cofinal sequence $(i_n:n<\omega)$ in
$\theta$ witnesses 
that a $(\mu,\theta)$-$d$-limit model is also a
$(\mu,\omega)$-$d$-limit model.%}\marginpar{PZ -abril 4 2015} 
\\ \\
\indent Therefore, we may assume for the rest of the proof that
$cf(\theta)>\omega$. 
\\ \\
\indent We construct by induction an $(\omega+1)$-sequence of
towers. We first construct by induction an $(\omega+1)$-sequence of
well-orderings $(I_n)_{n\leq \omega}$.
Let $I_0$ be a well-ordering with a sequence of
elements $(i_\alpha:\alpha\le \theta)$ in $I$ such that
$(i_\alpha:\alpha<\theta)$ is a cofinal sequence in $I_0\setminus
\{i_\theta\}$ (therefore $cf(I\setminus \{i_\theta\})=\theta$)
and $i_\theta$ is the maximum of $I_0$, as in
assumption~\ref{In}. Based on this $I_0$,
build inductively a sequence
$I_n$ of well-orderings just as in assumption~\ref{In}. Finally,
$I_\omega=\bigcup_{n<\omega}I_n$. All the stages $I_n$ (for
$n<\omega$) clearly have uncountable cofinality ($\theta$) below their
maximum element $i_\theta$. The final well-ordering $I_\omega$ also has
 uncountable cofinality below $i_\theta$.
\\ \\
\indent Let $\delta:\mathbb{R}^+\to \mathbb{R}^+$ be a mapping as in
definition~\ref{reduced_towers}. First consider any $I_0$-s-tower
$\tower^0$  such that
$M_0^0:=M$.  Suppose now that we have constructed
$\tower^n$
%\in \C{K}_{\mu,I_n}$
 an $I_n$-tower. By lemma \ref{weak_full} and proposition
\ref{reduced_density}, there exists an s-tower $\tower^{n+1}\ge
\tower^{n}$ which is $d$-reduced relative to $\delta$ and also
satisfies the properties given
in lemma \ref{weak_full}. At stage $\omega$, let
$\tower^{\omega}$ be
the completion of the union of $\langle \tower^n:n\le m
<\omega\rangle$; this is an $I_\omega$-tower. By proposition
\ref{Union_reduced},
$\tower^\omega$ is a d-reduced tower relative to the mapping
$\delta':\mathbb{R}^+\to \mathbb{R}^+$ defined at the end of the proof
of proposition~\ref{Union_reduced} (hence continuous, by
proposition~\ref{reduced_continuity}). Since
  $i_\theta=\max(I_0)$, then $i_\theta$ is also the $\max$ of each one of
  the $I_n$'s.
Notice that $M_{i_\theta}^\omega$ is a $(\mu,\omega)$-$d$-limit model
witnessed by $\{M_{i_\theta}^n:n<\omega\}$: this holds by the
definition of $\le$. Also, notice that $M_{i_\theta}^\omega$ is a
$(\mu,\theta)$-$d$-limit model. [Why?
$\tower^\omega$ is relatively full with respect to
$(M^n_i)_{n<\omega,i\in
  I_\omega}$ (by lemma \ref{weak_full}); so, by proposition
\ref{LimitModels}, $M_{i_\theta}^\omega$ is a $(\mu,\theta)$-$d$-limit
witnessed by $\{M_{i}^\omega:i<i_\theta \}$ (notice that continuity
of reduced towers guarantees that
$M_{i_\theta}^{\omega}=\overline{\bigcup_{i<i_\theta}
  M_{i}^\omega}$)].

So, we have constructed a $(\mu,\omega)$-$d$-limit model over $M$
which is also a $(\mu,\theta)$-$d$-limit model over $M$.
\edem[Prop. \ref{omega_theta}]

With this, we finish putting together the proof of our main theorem.

\begin{theorem}\label{Uniqueness}%\index{Materias}{d limit
                                %model@$d$-Limit Model!Uniqueness}
Let $\C{K}$ be an MAEC satisfying AP, JEP, CTP, existence of large
enough models, assumption~\ref{Superstability} (superstability) and
which is $\mu^+$-categorical.
If $M_i$ is a $(\mu,\theta_i)$-$d$-limit over $M$ ($i\in \{1,2\}$),
then $M_1\approx_M M_2$.
\end{theorem}
\bdem
By proposition \ref{omega_theta}, any two such models $M_1$ and $M_2$
are both isomorphic (over $M$) to an arbitrarily picked
$(\mu,\omega)$-$d$-limit model $M_0$ over $M$. They are therefore
isomorphic over $M$. 
\edem[Theorem \ref{Uniqueness}]

\end{document}